\documentclass[a4paper,11pt]{article}
\pdfoutput=1
\usepackage{graphicx}  %%% for including graphics
\usepackage{url}       %%% for including URLs
\usepackage{times}
\usepackage[square, numbers]{natbib}
\usepackage[margin=25mm]{geometry}

\usepackage{parskip}

\usepackage{multirow}
\usepackage{slashbox}
\usepackage{multicol}
\usepackage{tikz}

\usepackage{diagrams}
\usepackage{epsfig}
\usepackage{latexsym}

\usepackage{savesym}
\usepackage{amssymb}
\usepackage{amsthm}
\usepackage{amsmath}
\savesymbol{iint}
\usepackage{txfonts}
\restoresymbol{TXF}{iint}

\usepackage{amssymb}

\newtheorem{theorem}{Theorem}[section]

\theoremstyle{definition}
\theoremstyle{definition}
\theoremstyle{definition}\newtheorem{definition}[theorem]{Definition}
\theoremstyle{definition}
\theoremstyle{definition}
\theoremstyle{definition}
\theoremstyle{definition}
\theoremstyle{definition}
\theoremstyle{definition}

\newcommand{\bit}{\begin{itemize}}
\newcommand{\eit}{\end{itemize}\par\noindent}
\newcommand{\ben}{\begin{enumerate}}
\newcommand{\een}{\end{enumerate}\par\noindent}
\newcommand{\beq}{\begin{equation}}
\newcommand{\eeq}{\end{equation}\par\noindent}
\newcommand{\beqa}{\begin{eqnarray*}}
\newcommand{\eeqa}{\end{eqnarray*}\par\noindent}
\newcommand{\beqn}{\begin{eqnarray}}
\newcommand{\eeqn}{\end{eqnarray}\par\noindent}

\usepackage{color}  

\def\bR{\begin{color}{red}}
\def\bB{\begin{color}{blue}}
\def\bM{\begin{color}{magenta}} 
\def\bC{\begin{color}{cyan}}
\def\bW{\begin{color}{white}}
\def\bBl{\begin{color}{black}}
\def\bG{\begin{color}{green}} 
\def\bY{\begin{color}{yellow}}
\def\e{\end{color}}

\def\lto{\multimapinv}

\title{Lambek vs.~Lambek: Functorial Vector Space Semantics\\ 
and String Diagrams for Lambek Calculus}

\author{
{\bf Authors:} Bob Coecke, Edward Grefenstette, Mehrnoosh Sadrzadeh\\ \ \\
{\bf Affiliation:} {Department of Computer Science, University of Oxford}\\ \ \\
{\bf Address:} Wolfson Building, Parks Road, Oxford, OX1 3QD\\ \ \\
{\bf Emails:} {\tt firstname.lastname@cs.ox.ac.uk}
}   
\date{}

\begin{document}
\maketitle

\begin{abstract} 

The Distributional Compositional Categorical (DisCoCat)  model  is a mathematical framework that provides  compositional semantics for meanings of natural language sentences. %Notably, it 
It consists of a computational procedure for constructing  meanings of sentences, given their grammatical structure  in terms of compositional type-logic,  and  given the empirically derived meanings of their words.  For the particular case that  the meaning of words  is modelled within a  distributional vector space model,  its experimental predictions, derived from real large scale data,  have outperformed  other empirically validated methods that could build vectors for  a  full sentence.   This success can be attributed to a conceptually motivated mathematical underpinning, something which the other methods lack, by integrating \em qualitative \em compositional type-logic and \em quantitative \em modelling of meaning within a category-theoretic mathematical framework. 
% but which  moreover did not have conceptually motivated mathematical underpinnings.  
 %, inspired by the use in categorical quantum mechanics.
%is inspired by the categorical vector space models of quantum protocols and ~\cite{AbrCoe2}. 
%The framework is obtained by combining two orthogonal approaches to natural language semantics: the distributional vector space models that represent  empirical meanings of words and the compositional type-logics that model the grammatical structure of sentences. 

% This combination is achieved  by taking the cartesian product of the  category of a pregroup type-logic and finite dimensional vector spaces, which are both compact closed. The combination  can be recasted by interpreting the pregroup in the vector spaces via a functor, as done in~\cite{PrellerSadr}.  The  functorial passage to the category of vector spaces, is well known as  `quantisation' in Topological Quantum Field Theory (TQFT).  Hence by analogy, here we are quantizing the grammar and developing a `grammatical quantum field theory'. 
The type-logic  used in  the DisCoCat  model  is Lambek's pregroup grammar. Pregroup types form a posetal compact closed category, which can be passed,  in a functorial manner,  on to the compact closed structure 
%A pregroup shares a  compact closed structure with the category 
of vector spaces,  linear maps  and  tensor product. 
 The diagrammatic versions of the equational reasoning in compact closed categories 
%This shared  structure allows for reasoning about 
can be interpreted as the  \emph{flow of word meanings} within sentences. % in terms of the diagrammatic calculus of compact closed categories. 
Pregroups  simplify 
%are a simplification of  
Lambek's previous type-logic, the Lambek calculus. 
The latter and its extensions have been  extensively  used to formalise  and reason about various linguistic phenomena. Hence,  the apparent reliance of the DisCoCat on pregroups has 
%drawn criticism from the linguistics community and 
 been seen as  a 
%its 
shortcoming.  
This paper addresses this concern,  by pointing out that one may as well realise a functorial passage from  the original type-logic of Lambek, a   monoidal bi-closed  category,  to vector spaces, or to any other model of meaning organised within a monoidal bi-closed  category. 
%In this paper we show that 
%is also compatible with the formalism of a DisCoCat.  This is done by developing a functor from  the Lambek calculus to the category of finite dimensional vector spaces. 
The corresponding string diagram calculus, due to  
%We  show how the graphical calculus of 
Baez and Stay, now depicts the  flow of word meanings, 
 and  also reflects the structure of the parse  trees of the Lambek calculus. 
%based on the information encoded in their   parse trees.  

%But there is a caveat in the use of pregroups versus Lambek calculi in composing meaning: the latter are not adequate at the level of expressing  meaning. 

%: the monoidal structure of vector spaces does not provide us with the same vector meaning as its compact structure. 
%and that further more this  setting and calculus provide more information about the constitutent structure of language to our semantic model than was previously available.
	\end{abstract}

\newpage

\section{Introduction}
Language is both empirical and compositional: we learn  meanings of  words by being exposed to linguistic practice, and we form sentences we've never heard before by composing words  according to the rules of  grammar. Various mathematical and formal models have sought to capture facets and aspects of language learning and formation. Compositional type-logical approaches~\cite{Lambek} represent sentence formation rules based on formal syntactic analysis, using formalisms such as context free grammars~\cite{chomskyCFG,montaguebook}, Lambek calculus~\cite{Lambek}, or Combinatorial Categorial Grammar \cite{steedmanCCG}.  Such formal approaches to grammar align well with Frege's notion of compositionality, according to which the meaning of a sentence is a function of the meaning of its parts~\cite{fregeSR}, but eschew the empirical nature of language, requiring pre-defined mathematical structures, domains and valuations to make sense. 

Orthogonal to formal logical models, empirical approaches to semantics construct representations of individual words based on the contexts in which they are used. These models are often referred to as \emph{distributional} or geometric models of semantics and are sometimes considered to be  in line with the ``meaning is use'' view of Wittgenstein's philosophy of language \cite{wittgensteinPI}. Distributional models have been applied successfully to tasks such as thesaurus extraction \cite{Grefenstette,Curran}, automated essay marking \cite{landauerdumais1997}, and other semantically motivated natural language processing tasks.  While these models reflect the empirical aspects of language learning that type-logical models lack, they in turn  lack composition operations which would allow us to learn meanings of phrases based on  the meanings of their parts. Developing models that could combine the strengths of the above  two approaches has proved to be a challenge for computational linguistics and its applications to natural language processing (NLP).

The  distributional compositional categorical (DisCoCat) model of meaning,  developed in 
%previous work 
\cite{CCS,CSC}, provides a solution to the above problem.  This framework,  which realised the challenge proposed in \cite{ClarkPulman} enables  a combination of the  type-logical and distributional models of meaning  and resulted in a procedure for  compositionally computing   meaning vectors for sentences by exploiting the grammatical  structure of sentences and the meaning vectors of the words therein.  
The framework was inspired by the category-theoretic  high-level framework for modelling quantum protocols~\cite{AbrCoe2},  were  the corresponding string diagram calculus exposes flows of information between the systems involved in  multi-system protocols such as quantum teleportation \cite{Kindergarten}. 
%,  was a cover heading feature in New Scientist (11 Dec.~2011) \cite{NewScientist}, 
 The DisCoCat model  has meanwhile been  experimentally validated  
 for 
%via its application to  
natural language tasks such as word-sense disambiguation within phrases \cite{GS,GS2}. 
%\footnote{EMNLP is the leading conference on corpus-based experimental NLP.}.

The  DisCoCat  model  relies on  Lambek pregroups~\cite{Lambek6} as its base  type-logic. 
 In category-theoretic terms, these have a 
%convenient 
(non-symmetric) compact closed structure  when considering types as objects and type reductions as morphisms. The DisCoCat exploits the fact that finite dimensional vector spaces can also be organised within a compact closed category.  The first and mainly technical goal  of this paper is to stress that the choice of  a compact or monoidal 
type-logic is  not crucial  
%irrelevant 
to the applicability of the procedure.   To achieve this goal,  we have  tweaked 
%extended 
the distributional compositional model of previous work from Lambek pregroups to  Lambek monoids,   hence  developing a  vector space model for the meaning of natural language sentences parsed within the  Lambek calculus. In this paper, we develop a similar homographic passage via a functor from a monoidal bi-closed category of grammatical types and reductions to the symmetric monoidal closed category of finite dimensional vector spaces. 

This functorial passage  is another contribution of this paper and gives rise to an interesting analogy with  Topological Quantum Field Theory (TQFT)~\cite{Atiyah, BaezDolan, Kock}. A TQFT is also a monoidal functor from the category of cobordisms into the category of vector spaces and linear maps.  From the perspective of TQFT, our DisCoCat models  form a   `Grammatical Quantum Field Theory' obtained by replacing the  monoidal  category of cobordisms in a  TQFT  by a certain partially ordered monoid which accounts for grammatical structure. This analogy of the compositional  distributional model of meaning with TQFT was first pointed out by Louis Crane at a workshop in Oxford, August 2008. Similar to the original model-theoretic framework of  meaning  by Montague~\cite{montaguebook}, this  semantic framework is  obtained via a homomorphic  passage from sentence formation  rules to compositions of meanings of words. However, contrary to  the Montague's model,  meanings of words and sentences are expressed in terms of vectors and vector compositions rather than in terms of sets and set-theoretic operations.

As a result of the  compactness of Lambek pregroups,    the mechanism of how  meanings of words interact to produce  meanings of sentences has a purely diagrammatic form, which admits an intuitive interpretation in terms of \emph{information flow}.  By \emph{information flow} we  mean the topology of the two-dimensional graphical representation of the operations that produces the meaning of sentences from the meaning of words.   Mathematically, these are expressed in the graphical language of the particular  category in which we model the meaning of words and sentences \cite{SelingerSurvey}, a practice tracing back to Penrose's work in the early 1970s \cite{Penrose},  that was  turned into a formal discipline by Joyal and Street in the 1990s \cite{JoyalStreet}.  These diagrams, for the particular case of compact closed categories,  were extensively exploited in the earlier  DisCoCat models. 
%The compact closed form of these diagrams were applied to reason about quantum protocols~\cite{AbrCoe2} and  natural language meaning~\cite{CCS,CSC}.  In this paper, we   apply these ideas to the monoidal approach and 
 Here we  show how  the clasp-string calculus of Baez and Stay~\cite{Baez} can be used to provide diagrams for information flows that arise in the Lambek monoids, which are not  compact.    Our ambition is to use this work as a starting point for providing vector space meaning for  more expressive  natural language sentences  such as those parsed with Combinatorial Categorial Grammars (CCGs) or  Lambek-Grishin calculus~\cite{Moortgat09,Bernardi}. The expressive powers of these  grammars go beyond that of Lambek grammars, which are context free.

%More concretely, to pass from a grammar to a category, the rules which model the type reductions are now accompanied by a map which  propagates the meaning of the words/phrases of that type into the meaning of the whole sentence. %Rather than arrows/graph there is also a `wave of meaning' flowing along that graph, ... 
%The title of this paper is a dedication to Lambek's seminal contributions to the field of type grammars, some of which are half a century apart, and his pioneering contributions to categori\underline{c}al logic, an area which very much inspired the above mentioned algorithm.  

%Traditionally  these kinds of functorial passages are used in TQFT and are referred to as \emph{quantisation}. Hence, our work in providing distributional meaning given a type grammar can be seen as quantizing the syntax of natural language. 

%for Ajdukiewicz-Bar-Hillel grammars \cite{Ajdukiewicz,Bar-Hillel}, Grishin grammars \cite{Grishin} as well as the even weaker grammars modeled by protogroups \cite{Lambek6}, \bM and even non-associative generalizations of any of these, see e.g.~\cite{Moortgat}\e.  

%The 2nd point is an investigation in what the choice of a type grammar (and also, the way in which it is used) actually implies at the level of the information flows of word meanings in sentences (of course, besides the constraints it imposed on well-typedness of sentences). 

Finally, drawing a connection with games seems appropriate in the context of this special issue. Application of games to interpreting and formalising natural language traces back to  the `dialogical logic' of Lorenz and Lorenzen~\cite{lorenzen} who used the dialogue analogy to develop  a game semantic model for formulae of intuitionistic logic. Later, a classical logic version of the theory with a model theoretic focus was developed  by Hintikka and Sandu~\cite{hintikka}.  A proof-theoretic approach  led to the use of linear logic, and proof nets, e.g. see~\cite{Lamarche,lecomte}. Independently, another line of research was pursued by  linguists who also used the term `dialogue games' to provide a semantic model for  real-life human-computer dialogues. One of the original  proposals of this line  was based on Grice's pragmatic philosophy of language and used  component programs and specifications to model dialogues and queries; the setting was applied  to online sale tools~\cite{levin}. Later on, a formal  model based on belief revision and Bayesian update was  developed for this approach~\cite{pulman}.   Our work can be seen as bridging these two (abstract logical and applied linguistic) communities. Our starting point is a  Lambek calculus, with a proof theory similar to that of intuitionistic multiplicative linear logic. In this calculus,  the grammatical structure of a sentence is represented as a derivation in a proof tree and depicted in diagrams similar to  proof nets. We interpret these derivations in vector spaces, seen as a monoidal closed category, and depict the grammatical and semantic interactions  via Baez-Stay diagrams. General linear logic proof nets are quite different from Baez-Stay diagrams, but their compact variants introduced in~\cite{ross,abramsky} resemble the compact closed string diagrams used in the pregroup derivations.

\section{Partial Order Structures in  Linguistics}
Application of partially ordered algebras to linguistics originated in the seminal work of Lambek in the 50's~\cite{Lambek}. In his debut work, Lambek  showed how a partially ordered residuated monoid can be used to analyze the syntactic structure of a fragment of English. He later  developed a decision procedure for this setting, based on a cut-free sequent calculus. This calculus can be seen as the father of linear logic, as it shares with it a monoidal tensor, in a non-commutative form,  and hence  two linear implications.  Lambek's approach was based on  a partial order of grammatical types; similar ideas have been present in the work of Bar-Hillel~\cite{Bar-Hillel} but were not formalized in algebraic and proof theoretic forms. About half a century later  in the late 90s, Lambek simplified his original residuated monoids in favor of a simpler partial order, which he called a pregroup~\cite{LambekBook}. Pregroups have been applied to analyse  various different languages, for references see~\cite{CasadioLambek}. In this section we review these two structures and their application to natural language.

\subsection{Lambek Monoids}
Lambek calculus~\cite{Lambek} is usually a reference to Lambek's  sequent calculus; this calculus is similar to that of intuitionistic multiplicative bi-linear logic, but lacks negation. It has  one main binary operation, which is noncommutative,  hence  has a right and a left implication. This calculus is sound and complete with regard to  partially ordered residuated monoids.  

Recall that for two order-preserving maps $f:A\to B$ and $g: B\to A$ on two partially ordered sets $A$ and $B$,  we say that $f$ is the \em left adjoint \em to, or the left \emph{residual} of, $g$ (or equivalently, $g$ is the \em right adjoint \em to $f$ or its left \emph{residual}), denoted $f \dashv g $, iff 
\[
\forall a\in A, b\in B, \quad f(a)\leq b \Leftrightarrow a \leq g(b)
\]
The above  condition is  equivalent to the following:
\[
\forall b\in B,  \quad f(g(b))\leq b, \qquad \text{and} \qquad \forall a\in A,  \quad a\leq g(f(a))
\]
 Based on the above definition, a residuated  monoid is defined as follows:
\smallskip
\begin{definition}
A  residuated monoid $(L, \leq, \cdot, 1, \multimap, \lto)$ is a partially ordered set $(L, \leq)$, equipped with a  monoid structure $(L, \cdot, 1)$ that preserves the partial order, that is for all $a,b, c \in L$, we have:
\[
\text{If} \quad a \leq b \quad \text{then} \quad a \cdot c\leq b \cdot c \quad \text{and} \quad c \cdot a \leq c \cdot b
\]
The unit element 1 satisfies the following  for all $a \in L$: 
\[
1 \cdot a =  a \cdot 1 = a
\]
That the monoid is residuated means that $\multimap $ and $\lto$  are the two adjoints of $\cdot$,  that is, 
$a\cdot (-)\dashv a \multimap (-)$ and $(-)\cdot b\dashv (-) \lto b$, explicitly we have:
\beq
b\leq a \multimap c \Leftrightarrow a\cdot b \leq c \Leftrightarrow a \leq c\lto b\,,
\eeq
or, equivalently, using the corollaries of these adjunctions, we have:
\beq
a \cdot (a \multimap c)\leq c\quad \mbox{,} \quad c\leq  a \multimap (a \cdot c) \qquad \mbox{,} \qquad   
(c\lto b)\cdot b \leq c\quad , \quad c\leq   (c \cdot b)\lto b\ ,
\eeq
\end{definition}

These structures are also referred to as  \emph{residuated lattices} in the literature. But strictly speaking, a residuated monoid  is a residuated lattice  with the exclusion of  its lattice operations, hence the main operation of this structure is a  monoid multiplication. 

We refer to a partially ordered  residuated monoid as a \emph{Lambek monoid}.  These monoids are applied to the encoding of the grammatical structure of natural language, whereby elements of the algebra denote grammatical types, the monoid multiplication is the juxtaposition of these types, and its unit is the empty type. The right and left adjoints are used to denote function-types; these encode the  types of the words that have a relational role, for example adjectives, verbs, adverbs,  conjunctives, and relative pronouns. This application procedure is formalised via the following structures.

\smallskip
\begin{definition}
For  $\Sigma$ the set of words of a natural language and ${\cal B}$ a set of basic grammatical types,  a Lambek type-dictionary is a binary relation $D$, defined as
\[
D \subseteq \Sigma \times T({\cal B})
\]
where $T({\cal B})$ is the free Lambek monoid generated over ${\cal B}$ (for the free construction see~\cite{Lambek}). 
\end{definition}
\begin{definition}
A Lambek grammar $G$ is a pair $\langle D, S\rangle$, where $D$ is a Lambek type-dictionary and $S \subset {\cal B}$ is a set of designated types, containing types such as that of a declarative sentence $s$, and a question $q$.  
\end{definition}

A Lambek grammar  is used to define the grammatical sentences of a language as follows.

\smallskip
\begin{definition}\label{grammar}
A string of words $w_1 w_2 \cdots w_n$ each of them from  $\Sigma$ is said to be grammatical iff for $1 \leq i \leq n$,  there exists a $(w_i, t_i)$ in the dictionary $D$, such that for the designated type  $s$ of a sentence in $S$,  the following partial order holds in $T({\cal B})$:
\[
t_1 \cdot t_2 \cdot \cdots \cdot t_n \leq s
\]
\end{definition}

As an example,  consider a simple language that contains the following words: 
\[
\Sigma = \left \{\text{men}, \text{dogs}, \text{cute}, \text{kill}, \text{do}, \text{not}\right\}
\]
For the sake of this example suppose that  in this language one can only form   declarative sentences `men kill dogs',  `men do not kill dogs',  and `men kill cute dogs'.  A Lambek type-dictionary that generates the grammatical sentences of this language has the following basic types: 
\[
{\cal B} = \{n,s, j, \sigma\}
\]
Here,  $n$ stands for a noun phrase,  $s$ for a declarative sentence, $j$ for the infinitive of a verb, and $\sigma$ is a marker type. The type assignments to words of $\Sigma$ are presented in Table~\ref{lambekdict}.

\begin{table}[h]
\begin{center}  
\begin{tabular}{c|c|c|c|c|c|c}
men & dogs & cute & kill  & to kill &  do &not \\
\hline
$n$ & $n$ & $n \lto n$ & $(n\multimap s) \lto n$  & $(\sigma \multimap j)\lto n$ & $(n \multimap s) \lto (\sigma \multimap j)$ &  $(\sigma \multimap j) \lto (\sigma \multimap j)$
\end{tabular}
\end{center}
\label{lambekdict}
\caption{Type Assignments for the Toy Language $\Sigma$ in a Lambek monoid.}
\end{table}

The Lambek type-dictionary $D$ corresponding to the type assignments of Table~\ref{lambekdict} is as follows:
\begin{align*}
&\Big\{(\text{men},n), (\text{dogs},n), (\text{cute}, n \lto n), (\text{kill}, (n\multimap s) \lto n), 
(\text{to  kill}, (\sigma \multimap j)\lto n), \\
\quad & \quad (\text{do}, (n \multimap s) \lto (\sigma \multimap j)), (\text{not}, (\sigma \multimap j) \lto (\sigma \multimap j)) \Big\}
\end{align*}
Note that the verb `kill' has two types, represented by two pairs in the type dictionary: the first one $(\text{kill}, (n\multimap s) \lto n)$ is for its transitive role, e.g.~in the sentence `men kill dogs' and  the second one $(\text{kill}, (\sigma \multimap j)\lto n)$ for its infinitive role, e.g.~in the sentence `men do not kill dogs'. 

By definition~\ref{grammar}, the sentence `men kill dogs' is  grammatical, since if we apply the monoid multiplication to  the types that correspond to the words, we obtain the term $n \cdot (n\multimap s) \lto n \cdot n$, which has the following reduction: 
\begin{align*}
n \cdot ((n\multimap s) \lto n) \cdot n &\leq \\
n \cdot (n\multimap s) &\leq s
\end{align*}
The sentence `men kill cute dogs' is also grammatical; it has the following reduction:
\begin{align*}
n \cdot ((n\multimap s) \lto n) \cdot (n \lto n) \cdot n &\leq\\
n \cdot ((n\multimap s) \lto n) \cdot  n&\leq\\
n \cdot (n\multimap s) &\leq s
\end{align*}
Similarly, the sentence `men do not kill dogs' is grammatical, according to the following reduction: 
\begin{align*}
n \cdot ((n \multimap s) \lto (\sigma \multimap j)) \cdot ((\sigma \multimap j) \lto (\sigma \multimap j)) \cdot ((\sigma \multimap j)\lto n) \cdot n &\leq\\
n \cdot ((n \multimap s) \lto (\sigma \multimap j)) \cdot ((\sigma \multimap j) \lto (\sigma \multimap j)) \cdot (\sigma \multimap j) &\leq\\
n \cdot ((n \multimap s) \lto (\sigma \multimap j)) \cdot (\sigma \multimap j) &\leq\\
n \cdot (n \multimap s)&\leq s\\
\end{align*}

\subsubsection{Lambek Pregroups}
In 1999, Lambek revisited his monoidal structures and introduced a simplification~\cite{Lambek6}. Instead of working in a partially ordered residuated monoid, he argued for a pregroup, which has  one non-residuated binary operation, but where each element of the partial order is required to have a left and a right adjoint. More precisely, we have:

\smallskip
\begin{definition} A Lambek pregroup is a partially ordered unital monoid where each element has a left and a right adjoint $(P, \leq, \cdot, 1, (-)^l, (-)^r)$. That is, for every $p \in P$, there is a $p^r$ and a $p^l$ in $P$, which satisfy the following four inequalities:
\begin{align*}
p\cdot p^r &\leq 1 \leq p^r \cdot p\\
p^l \cdot p &\leq 1 \leq p \cdot p^l
\end{align*}
\end{definition}
From this definition it follows that adjoints are unique and reverse the order, that is if we have $p \leq q$ for $p,q \in P$ then it follows that $q^l \leq p^l$ and also that $q^r \leq p^r$. One can also show that the unit is self adjoint, that is $1^r = 1 = 1^l$, that opposite adjoints cancel out, that is $(p^r)^l = (p^l)^r = p$, but same adjoints iterate, for instance $(p^r)^r$ is not necessarily equal to $p$ and neither is $(p^l)^l$. Another nice property is that the monoid multiplication is self adjoint, that is $(p \cdot q)^r = q^r \cdot p^r$ and also $(p \cdot q)^l = q^l \cdot p^l$. These properties are all proved and elaborated on by Lambek~\cite{Lambek6}. 

Apart from linguistics,  pregroups have concrete applications in other fields such as  to number theory~\cite{Lambek7}. Notably,  an example of a pregroup structure on natural numbers  is the set of all unbounded monotone functions on the set of integers $\mathbb{Z}$. Here, the partial order is the natural ordering of integers extended to functions, that is for  $f,g \in \mathbb{Z}^{\mathbb{Z}}$, we have:
\[
f \leq g  \quad \text{iff} \quad f(n) \leq g(n), \forall n \in \mathbb{Z}
\]
The monoid multiplication is the composition of functions and its unit is the identity function, that is for $n \in \mathbb{Z}$ we have:
\[
(f \cdot g)(n) = f(g(n)) \qquad \text{and} \qquad 1(n) = n
\]
The left and right adjoints of each function are computed by using their canonical definitions and  via the supremum and infimum operations on integers, again extended to functions, as follows:
\[
f^r(n) = \vee\{m \in \mathbb{Z}\mid f(m) \leq n\} \hspace{2cm}
f^l(n) = \wedge \{m \in \mathbb{Z}\mid n \leq f(m)\}
\]
As an example, take $f(x) = 2x$, then compute $f^r(x) = \lfloor x/2 \rfloor$ and $f^l(x) = \lfloor(x+1)/2 \rfloor$, where $\lfloor x \rfloor$ is the biggest integer less than or equal to $x$ (for details of these computations and more examples see~\cite{Lambek7}).

Pregroups are applied to the analysis of syntax in the same way as Lambek monoids. The types are translated from their monoidal implicative form to a pregroup adjoint form as follows:
\[
p \multimap q \quad \leadsto \quad p^r \cdot q \hspace{2cm} p \lto q \quad \leadsto \quad p \cdot q^l
\]
The pregroup version of  Table~\ref{lambekdict} is as follows:

\begin{table}[h]
\begin{center}  
\begin{tabular}{c|c|c|c|c|c|c}
men & dogs & cute  &  kill & to kill &do & not \\
\hline
$n$ & $n$ & $n \cdot n^l$ & $n^r\cdot s \cdot n^l$  & $\sigma^r \cdot j \cdot n^l$&$n^r \cdot s \cdot j^l \cdot \sigma$ &  $\sigma^r\cdot j \cdot j^l \cdot \sigma$
\end{tabular}
\end{center}
\label{pregdict}
\caption{Type Assignments for the Toy Language $\Sigma$ in a Lambek pregroup.}
\end{table}

The reduction corresponding to the sentence `men kill dogs' is as follows:
\[
n \cdot n^r \cdot s \cdot n^l \leq 1 \cdot s \cdot 1 = s
\]
For `men kill cute dogs', we have:
\[
n \cdot n^r \cdot s \cdot n^l \cdot n \cdot n^l \cdot n \leq 1 \cdot s \cdot 1 \cdot 1 = s
\]
And  for `men do not kill dogs' we have: 
\begin{align*}
n \cdot n^r \cdot s \cdot j^l \cdot \sigma \cdot \sigma^r\cdot j \cdot j^l \cdot \sigma \cdot \sigma^r \cdot j \cdot n^l \cdot n \leq
&\\ 
1 \cdot s \cdot j^l \cdot 1 \cdot j \cdot j^l \cdot 1 \cdot j \cdot 1 = &\\
 s \cdot j^l  \cdot j \cdot j^l\cdot j \leq&\\
s \cdot 1 \cdot 1 = s
\end{align*}
The pregroup reductions are usually depicted in \emph{cancelation diagrams} using  curved strings (cups), which connect the types that are being cancelled in each step, and lines, which depict  types which have not been canceled. Nested strings are used to denote multiple steps. For example, the cancelation diagram of `men kill dogs' is as follows ( the $\cdot$'s are dropped). 

\vspace{-2mm}
\begin{center}
\begin{minipage}{3cm} 
\hspace{-2cm}\begin{picture}(50,50)(150,150)
\put(200,185){men}
\put(208,170){$n$}
\put(230,185){kill}
\put(228,170){$n^r$}
\put(240,170){$s$}
\put(248,170){$n^l$}
\put(260,185){dogs.}
\put(270,170){$n$}
\end{picture}

\vspace{-0.5cm}
\hspace{-0.15cm}{\epsfig{figure=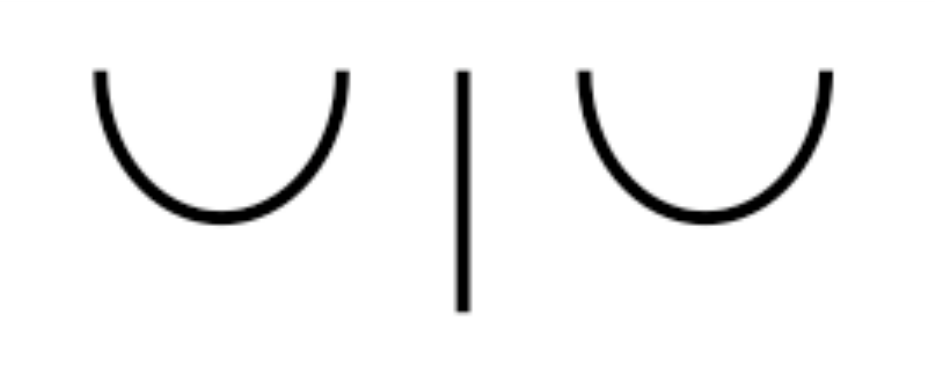,width=80pt}} 
\end{minipage}
\end{center}

The diagram for the grammatical structure of  `men do not kill dogs'  is as follows:

\vspace{0,5cm}
\begin{minipage}{6cm}
\begin{picture}(15,15)(150,148)

\put(292,165){men}
\put(300,150){$n$}

\put(340,165){do}
\put(330,150){$n^r$}
\put(338,150){$s$}
\put(343,150){$j^l$}
\put(351,150){$\sigma$}

\put(390,165){not}
\put(380,150){$\sigma^r$}
\put(389,150){$j$}
\put(395,150){$j^l$}
\put(402,150){$\sigma$}

\put(430,165){kill}
\put(425,150){$\sigma^r$}
\put(435,150){$j$}
\put(443,150){$n^l$}

\put(465,165){dogs.}
\put(470,150){$n$}
\end{picture}

\vspace{0.4cm}\hspace{5cm}
{\epsfig{figure=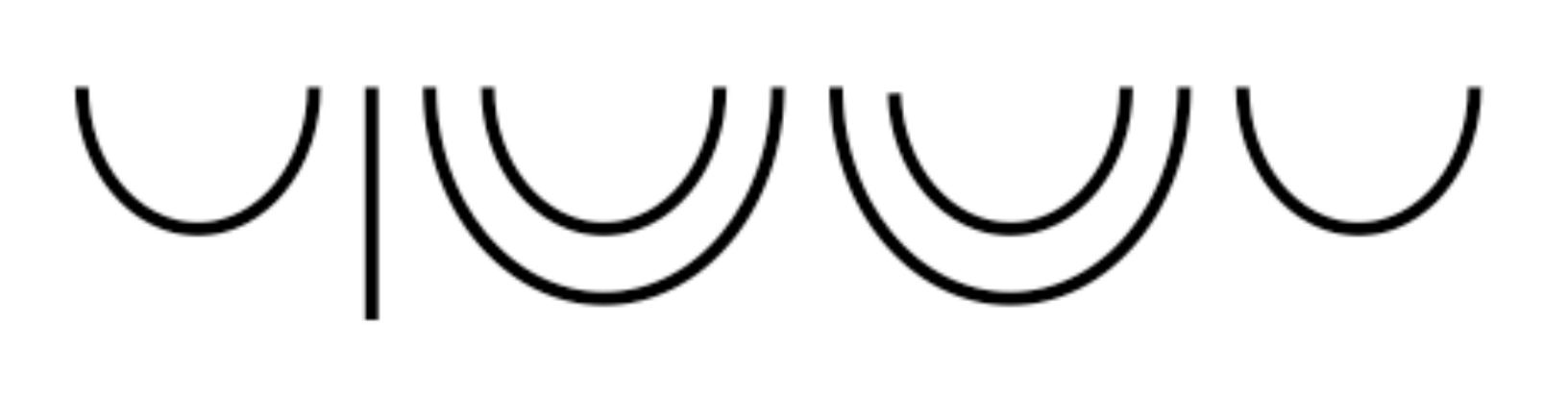,width=187pt}}
\end{minipage}

\noindent
These diagrams provide  an intuitive reading of the grammatical structure  of the sentence. For instance,  in the first diagram we read  that for a transitive sentence to be grammatical, the  verb has to interact with its subject and object, depicted via  the curved strings, hence producing   a sentence via the  line.

%Sequent calculi for pregroups have  been developed~\cite{Busz2}; proof nets for these are degenerate versions of the proof nets for compact bi-Linear Logic where tensor and par are the same. These are developed in the context of quantum protocols  and knot theory in~\cite{ross,abramsky}.  

The challenge of using type-logics  to formalize natural language grammar resides in assigning the right basic and compound types to the words of the language such that they would generate all the grammatical statements of the language and do not over-generate. Different type dictionaries have been suggested for different languages, e.g.~see~\cite{LambekBook,Moortgat}. The types used in this paper are from~\cite{Preller}; these are  chosen to keep the setting simple and to be able to parse our example sentences. Parsing a more complex language needs a  more elaborate type-dictionary. 

Additional operators, either in terms of modalities or new additive binary connectives have been added to type-logics  to increase their expressive power~\cite{MootThesis,Morrill,Moortgat96,kislak}. The expressive power of Lambek monoids and pregroups is the same as that of \emph{context-free} grammars of the Chomsky hierarchy~\cite{pentus,Busz}. The extensions with modalities and additives increase it to weak variants of \emph{context-sensitive} such as \emph{mildly context-sensitive}.  In the proof-theoretic calculi of these algebras,  the grammatical reductions are depicted via proof nets~\cite{Roorda,Moot}, which offer a richer analysis of the decomposition of types.

The above  type-logical structures do represent the  grammatical structure of a language in a compositional way, but   do not offer a model for  the   lexical semantics of the words of a language. In such type-logics, words are only modelled  by their grammatical roles, hence words that have the same grammatical role but different lexical meanings cannot be distinguished from one another. For instance,   the words `dog' and `men' are both  represented only by their grammatical role $n$,   ignoring the fact that they have  different lexical meanings. The same problem exists for transitive verbs such as `kill' and `eat',  adjectives such as `cute' and `green',  and so on.   For a subspace of a real vector space, see Figure~\ref{fig:VectMod}.

\section{Distributional Models of Meaning}

Orthogonal to the  type-logical models, distributional models of meaning are mainly concerned with lexical semantics. Best described by a quotation by Firth that ``You shall  know a word by the company it keeps."~\cite{Firth}, these   models are based on the dictum that words that often appear in the same context have similar meanings.  For example, the words `cat' and `dog' often appear in the context of  `pet', `furry', `owner', and `food', hence they have similar meanings. Similarly,  `kill' and `murder' often appear in the context of `police',  `gun', and `arrest',  hence these also have similar meanings.  

To formalize this idea, one builds a finite dimensional vector space $N$ whose basis vectors  are the context words. Ideally,  the context words are  all the lemmatized words of a corpus of interest.  Typically, these are reduced to a more refined set,  based on the domain of application.  One then fixes a window of $k$ words  and builds a vector for each word, representing its  lexical meaning. We denote the meaning  vector of a word  by   $\overrightarrow{\text{word}} = \sum_i {c_i \,\overrightarrow{n_i}}$, for $c_i$ a real number and $\overrightarrow{n}_i$ a basis vector of $N$. The $c_i$ weights are obtained by  first counting  how many times a word has  appeared within $k$ (e.g. 5) words of each context and then normalising this count.  The most popular normalisation measure is  Term-Frequency (TF) divided by Inverse-Document-Frequency (IDF); it assigns a degree of  importance to the appearance of a word in a document.   TF/IDF  is a proportion for the number of times a word appears in the corpus to the frequency of the total number of words in the  corpus.

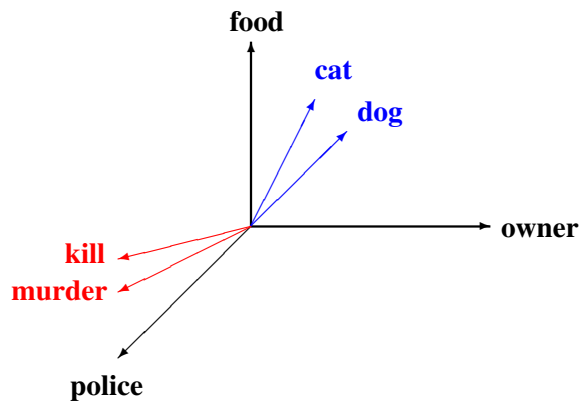
\begin{figure}[h]
\begin{center}
\begin{minipage}{10cm}
\begin{center}
\setlength{\unitlength}{0.7mm}
\begin{picture}(60, 70)
  \linethickness{0.3mm}
  \put(26,57){\text{\bf food}}
  \put(30, 20){\vector(1, 0){45}}
  \put(77,18){\text{\bf owner}}
  \put(30, 20){\vector(0, 1){35}}
  \put(30, 20){\vector(-1, -1){25}}
  \put(-4,-12){\text{\bf police}}
   
  \linethickness{0.3mm}
   
  \put(30, 20){\begin{color}{blue}\vector(1, 1){18}\end{color}}
  \put(42,48){\text{\bf \begin{color}{blue}cat\end{color}}}
  
  \put(30, 20){\begin{color}{blue}\vector(1, 2){12}\end{color}}
  \put(50,40){\text{\bf \begin{color}{blue}dog\end{color}}}
  
  \put(30, 20){\begin{color}{red}\vector(-2, -1){25}\end{color}}
  \put(-5,13){\text{\bf \begin{color}{red}kill\end{color}}}
  
  \put(30, 20){\begin{color}{red}\vector(-4, -1){25}\end{color}}
  \put(-15,6){\text{\bf \begin{color}{red}murder\end{color}}}
%  \put(30, 20){\begin{color}{magenta}\vector(-1, 2){14}\end{color}}
%  \put(-17,50){\text{\bf \begin{color}{magenta}gaze $\otimes$ gesture\end{color}}}
 \end{picture}
\end{center}

\vspace{0.9cm}
\end{minipage}
\end{center}
\caption{A Subspace  of a Real Vector Space Model of Meaning}
\label{fig:VectMod}
\end{figure}

The distance between the meaning vectors, for instance the cosine of their angle, provides a good measure of similarity of meaning. For example, in the vector space of Figure~\ref{fig:VectMod},  the angle between meaning vectors of `cat' and `dog'  is   small  and so is the angle between meaning vectors of  `kill' and `murder'. Various similarity measures have been implemented on large scale data (up to a billion words) to build high dimensional vector spaces  (tens of thousands of basis vectors). These have been successfully  applied to automatic generation of thesauri; for example  see~\cite{Curran}.

Contrary to the type-logical models, the empirical distributional models  ignore the grammatical structure of language. For instance, these models do not offer a canonical way of building  meaning vectors for  sentences. A compositional solution to this problem  should use  vector composition operators to combine  meaning vectors of words and obtain  meaning vectors for sentences. The simple operations widely studied in the literature, for instance in~\cite{Lapata},  are  addition ($+$) and component-wise multiplication ($\odot$). These are commutative and hence do not even respect the word order. If $\overrightarrow{vw} = \overrightarrow{v} + \overrightarrow{w}$ or $\overrightarrow{v} \odot \overrightarrow{w}$, then $\overrightarrow{vw} = \overrightarrow{wv}$, leading to unwelcome equalities such as  the following:
\[
\overrightarrow{\textrm{men kill dogs}} \ = \ \overrightarrow{\textrm{dogs kill men}}
\]
Inspired by the connectionist model of meaning in Cognitive Science~\cite{Smolensky}, a combination of  Kronecker product (which is non-commutative)  and syntactic relations has been suggested in~\cite{ClarkPulman} as a possible solution. The problem with this model is that the dimensionalities of sentence vectors  differ for sentences with different grammatical structures, barring them from being compared. For instance, in this model one cannot compare meanings of the two sentences `${\textrm{men kill dogs}}$' and `${\textrm{men kill}}$', since they live in different spaces.

\section{Distributional Compositional Categorical Model of Meaning}

In previous work~\cite{CCS,GS},  we combined the pregroup type-logic with the distributional model and developed a framework that produces vectors for meanings of sentences,  from their grammatical structure and the vectors of the words therein. This framework, summarised below,  is the Cartesian product of two compact closed categories: that of a pregroups with that of finite dimensional vector spaces.

\subsection{Pregroups and Vector Spaces as Compact Closed Categories}

%Recall the definition of a compact closed category. 

%\smalls kip
%\begin{definition}
Lambek's type-logics  closely relate to abstract categorical structures~\cite{Lambek89,PrellerLambek}.   To see this connection, we recall some definitions from category theory. A compact closed category has objects $A, B$, morphisms $f \colon A \to B$, a monoidal tensor  $A \otimes B$ that has a unit $I$, and  for each object $A$  two  objects $A^r$ and $A^l$ together with the following morphisms:
\begin{align*}
A \otimes A^r \rTo^{\epsilon^r} &I \rTo^{\eta^r} A^r \otimes A\\
A^l \otimes A \rTo^{\epsilon^l} &I \rTo^{\eta^l} A \otimes A^l\
\end{align*}
The above satisfy the following equalities, where $1_A$ is the identity morphism on object $A$:
\begin{align*}
(1_A \otimes \epsilon^l) \circ (\eta^l \otimes 1_A) & = 1_A \hspace{2cm} (\epsilon^r \otimes 1_A) \circ (1_A \otimes \eta^r) = 1_A\\
(\epsilon^l \otimes 1_A) \circ (1_{A^l} \otimes \eta^l) &= 1_{A^l} \hspace{2cm} (1_{A^r} \otimes \epsilon^r) \circ (\eta^r \otimes 1_{A^r}) = 1_{A^r}
\end{align*}
These inequalities are known as the \emph{yanking} equalities. Note  we do not assume symmetry of the tensor.
%\end{definition} 

\medskip
A pregroup is a compact closed category~\cite{PrellerLambek}, to which we refer as $Preg$. The elements of the partially ordered set $p, q \in P$ are the objects of $Preg$, the partial order relation provides the morphisms, that is,  there exists a morphism of type $p \to q$ iff $p \leq q$; we denote this morphism by $[p \leq q]$. Compositions of morphisms are given by transitivity and the identities by  reflexivity of the partial order. The monoid multiplication and  its unit (denoted by 1 rather than $I$) provide the tensor of the category, and the four adjoint inequalities provide the epsilon and eta morphisms, that is we have:
\begin{align*}
\epsilon^r = [p \cdot p^r \leq 1] &\hspace{2cm} \epsilon^l = [p^l \cdot p \leq 1]\\
\eta^r = [1 \leq p^r \cdot p] &\hspace{2cm}  \eta^l = [1 \leq p \cdot p^l]
\end{align*}

Finite dimensional vector spaces also form a compact closed category~\cite{AbrCoe2}, we refer to it as $FVect$.  Finite dimensional vector spaces $V, W$ are  objects of this category;  linear maps $f \colon V \to W$ are its morphisms, and composition is the composition of linear maps. The tensor product between the spaces $V \otimes W$ is the monoidal tensor whose unit is a field, in our case $\mathbb{R}$. As opposed to the tensor of the pregroup, this tensor is symmetric, hence we have a naturual isomorphism $V \otimes W \cong  W \otimes V$. As a result of symmetry of the tensor,  the two adjoints collapse to one and we  obtain  $V^l \cong V^r \cong V^*$, where $V^*$ is  the dual of  $V$. Since the basis vectors of our vector spaces are fixed, we furthermore obtain that $V^* \cong V$. Finally, a vector $\overrightarrow{v} \in V$ is represented by the morphism $\mathbb{R} \rTo^{\overrightarrow{v}} V$. 

Given a basis $\{r_i\}_i$ for a vector spaces $V$, the epsilon maps are given by the  inner product extended by linearity, that is we have:
\begin{align*}
\epsilon^l  =  \epsilon^r \colon \quad  &V^* \otimes V \to \mathbb{R}\\
:: \quad &\sum_{ij} c_{ij} \ \psi_i \otimes \phi_j  \quad \mapsto \quad \sum_{ij} c_{ij} \langle \psi_i \mid \phi_j \rangle
\end{align*}
Eta maps are %create what is referred to as \emph{maximally entangled states} or Bell states, 
defined as follows:
\begin{align*}
\eta^l = \eta^r \colon \quad & \mathbb{R} \to V \otimes V^*\\
::\quad  &  1 \quad \mapsto \quad \sum_i r_i \otimes r_i
\end{align*}
Here $1 \in \mathbb{R}$ is the number 1, and the above assignment extends to all other numbers by linearity. 

A DiscoCat is the cartesian product of $FVect$ and $Preg$. This is a category: its objects are pairs $(W, p)$ for $W$ a finite dimensional vector space and $p$ a pregroup type. Its morphisms are pairs of morphisms $(f \colon V \to W, [p \leq q])$, for $f$ a linear map and $p \leq q$ a pregroup partial order. Composition is obtained in a pointwise fashion by composing  linear maps and transitivity of the order. The identity morphism for an object $(V,p)$ is  a pair of identity morphisms $(1_V, [p\leq p])$ from $FVect$ and $Preg$. 

\medskip
\begin{definition}
The distributional compositional categorical (DisCoCat) model of meaning is the  category $FVect \times Preg$, equipped with a tensor product given by point wise tensor of $Preg$ and $FVect$, that is $(V, p) \otimes (W, q) = (V \otimes W, p \cdot q)$, whose unit is $(\mathbb{R}, 1)$. 
\end{definition}
It has been shown in~\cite{CCS} that the compact structure carries over from $Preg$ and $FVect$ to a DisCoCat, with epsilon and eta maps given in a  pointwise fashion as follows:
\begin{align*}
\epsilon^r = (V \otimes V^* \to \mathbb{R}, [p \cdot p^r \leq 1]) &\hspace{2cm} \epsilon^l = (V^* \otimes V \to \mathbb{R}, [p^l \cdot p \leq 1])\\
\eta^r = (\mathbb{R} \to V^* \otimes V, [1 \leq p^r \cdot p]) & \hspace{2cm} \eta^l = (\mathbb{R} \to V \otimes V^*, [1\leq p \cdot p^l])
\end{align*}

\subsection{Meaning Vectors  for Strings of Words}

The pair  $(\overrightarrow{w} \in W,p)$ is an object of $Fvect \times Preg$. We use this object to represent  the \emph{semantics} of a word $w$ and refer to it as the \emph{meaning} of the word $w$.  It consists of a vector space $W$ where the  meaning vector $\overrightarrow{w}$ of word  $w$ that  has  pregroup type $p$ lives.  Based on this notion, the  meaning vector representing the semantics of a sentence is obtained according to the following definition.

\medskip
\begin{definition}\label{meaningdef}\em
We define the  meaning vector \ $\overrightarrow{w_1 \cdots w_n}$ \ of  a string of words \ $w_1 \cdots w_n$ \ to be: 
\[
\overrightarrow{w_1 \cdots w_n} := f(\overrightarrow{w_1} \otimes \cdots \otimes \overrightarrow{w_n})
\]
where for  $(\overrightarrow{w_i} \in W_i, p_i)$  meaning   of  the word $w_i$,  the linear map $f$ 
is built by substituting each $p_i$ in  the pregroup reduction map of the string  \hbox{$[p_1 \cdot \cdots \cdot p_n \leq x]$} with  $W_i$. 

Thus for $\alpha = [p_1 \cdot  \cdots \cdot p_n \leq x]$ a morphism in $Preg$ and $f = \alpha[p_i/ W_i]$ a linear map in $FVect$, the following  is a morphism in  $FVect\times Preg$: 
\[
(W_1 \otimes \cdots \otimes W_n, \,p_1 \cdot \cdots \cdot p_n) \rTo^{(f, \leq)} (X,x)
\]
\end{definition}

For  example, to assign a meaning vector to an adjective-noun phrase, we start with the meanings of adjective and noun, which have the following forms:
\[
(\overrightarrow{\text{adj}} \in W \otimes W, n \cdot n^l) \hspace{2cm}
(\overrightarrow{\text{noun}} \in W, n)
\]
Here, we are assigning the vector space $W$ to the noun $\overrightarrow{\text{noun}}$ and assuming that the distributional meaning vectors of the nouns live in it, that is $\overrightarrow{\text{noun}} \in W$.  For  the  meaning vector  of the adjective, it is to assumed to be an element of $W \otimes W$, hence representable by $\sum_{lm} c_{lm}\ \overrightarrow{w}_l \otimes \overrightarrow{w}_m$,  for $\overrightarrow{w}_l, \overrightarrow{w}_m$ basis vectors of $W$.    The  meaning vectors of  the noun and adjective are represented by the following morphisms:
\[
\mathbb{R} \rTo^{\overrightarrow{\text{noun}}} W \hspace{2cm}
\mathbb{R} \rTo^{\overrightarrow{\text{adj}}} W \otimes W
\]

The pregroup reduction map of the adjective-noun phrase is as follows:
\[
\alpha = [n \cdot n^l \cdot n \leq n] = 1_n \cdot \epsilon_n
\]
Substituting each type in $\alpha$ with the  vector space associated with the words of that type, we obtain  the linear map $f$ corresponding to $\alpha$ to be the morphism  $(1_W \otimes \epsilon_W)$. 
The meaning vector of an adjective-noun phrase is computed by applying the linear map corresponding to the pregroup parse of the adjective-noun phrase, that is  $(1_W \otimes \epsilon_W)$ , to the tensor product of the meaning vector of the adjective $\overrightarrow{\text{adj}}$ with the meaning vector of the noun $\overrightarrow{\text{noun}}$. The  categorical  morphism corresponding to  this computation is as  follows: 
\begin{eqnarray*}
\overrightarrow{\text{adj \ noun}} & :=&  (1_W \otimes \epsilon_W) (\overrightarrow{\text{adj}} \otimes \overrightarrow{\text{noun}})\\
&\cong& (1_W \otimes \epsilon_W) \left(\left( \mathbb{R} \rTo^{\overrightarrow{\text{adj}}} W \otimes W\right) \otimes \left (\mathbb{R} \rTo^{\overrightarrow{\text{noun}}} W\right)\right)\\
&\cong& \mathbb{R} \otimes \mathbb{R} \rTo^{\overrightarrow{\text{adj}} \otimes \overrightarrow{\text{noun}}} W \otimes W  \otimes W  \rTo^{1_W \otimes \epsilon_W}W \otimes \mathbb{R} \\
&\cong& \mathbb{R}  \rTo^{\overrightarrow{\text{adj}} \otimes \overrightarrow{\text{noun}}} W\otimes W   \otimes W\rTo^{1_W \otimes \epsilon_W}W  
\end{eqnarray*}
The concrete value of the vector corresponding to the above morphism is:
\[
\overrightarrow{\text{adj \ noun}}  \quad:= \quad 
(1_W \otimes \epsilon_W) \left(\left(\sum_{lm} c_{lm}\ \overrightarrow{w}_l \otimes \overrightarrow{w}_m\right) \otimes \overrightarrow{\text{noun}}\right) \quad = \quad
\sum_{lm} c_{lm}\,  \overrightarrow{w}_l \, \langle \overline{w}_m \mid \overrightarrow{\text{noun}}\rangle\,
\]

The meaning vector of a  sentence with an adjective-noun phrase  is computed  by substituting the above for the meaning vector of the adjective-noun phrase when computing  the meaning of the sentence. For instance, consider the sentence `men kill dogs'; we have the following for the   meaning vectors of  words: 
\[
(\overrightarrow{\text{men}} \in W, n), \quad (\overrightarrow{\text{dogs}} \in W, n), \quad (\overrightarrow{\text{kill}} \in W\otimes S \otimes W, n^r \cdot s \cdot n^l)
\]
So we are  assigning the vector space $W$ to nouns and the vector space $S$ to  sentences; as a result   the vector space of  transitive verbs and in particular of `kill' will become  $W \otimes S \otimes W$; each verb being  representable by   $\sum_{ijk} c_{ijk} \overrightarrow{w}_i \otimes \overrightarrow{s}_j \otimes \overrightarrow{w}_k$, for $\overrightarrow{w}_i, \overrightarrow{w}_k$  basis vectors of $W$ and  $\overrightarrow{s_j}$ a basis vector of $S$.   The pregroup reduction map of the sentence is as follows:
\[
\alpha = [n \cdot n^r \cdot s \cdot n^l \cdot n \leq s] = \epsilon_n \cdot 1_s \cdot \epsilon_n
\]
Substituting each type in $\alpha$ with the  vector space associated with the words of that type, we obtain  $f$ to be the morphism $\epsilon_{W} \otimes 1_{S} \otimes \epsilon_W$. Given the distributional meanings of each word, that is  $\overrightarrow{\text{men}}, \overrightarrow{\text{kill}}, \overrightarrow{\text{dogs}}$,  the meaning of the sentence will be:
 \[
\overrightarrow{\mbox{men kill dogs}}  =  (\epsilon_{W} \otimes 1_{S} \otimes \epsilon_W) \left(\overrightarrow{\text{men}} \otimes \overrightarrow{kill} \otimes \overrightarrow{dogs}\right) \colon 
 W \otimes W \otimes S \otimes W \otimes W \rTo^{\epsilon_W \otimes 1_S \otimes \epsilon_W} S
 \]
further simplified as follows: 
 \begin{align*}
\overrightarrow{\mbox{men kill dogs}} =\ & f(\overrightarrow{\text{men}} \otimes \overrightarrow{kill} \otimes \overrightarrow{dogs}) \\
=\ &(\epsilon_W \otimes 1_S \otimes \epsilon_W) \left(\overrightarrow{\text{men}} \otimes\overrightarrow{\text{kill}} \otimes \overrightarrow{\text{dogs}}\right)\\
=\ & (\epsilon_W \otimes 1_S \otimes \epsilon_W) \left(\overrightarrow{\text{men}} \otimes \left(\sum_{ijk} c_{ijk} \overrightarrow{w}_i \otimes \overrightarrow{s}_j \otimes \overrightarrow{w}_k \right) \otimes \overrightarrow{\text{dogs}}\right)\\
=\ & \sum_{ijk} c_{ijk} \, \langle \overrightarrow{\text{men}} \mid \overrightarrow{w}_i \rangle  \langle \overrightarrow{w}_k \mid \overrightarrow{\text{dogs}} \rangle \ \overrightarrow{s}_j 
\end{align*}
Now, we can compute the meaning vector of the sentence `men kill cute dogs':  the pregroup morphism of the grammatical reduction of this sentence is $\epsilon_n \cdot 1_s \cdot \epsilon_n \cdot \epsilon_n$, whose linear map will be $\epsilon_W \otimes 1_S \otimes \epsilon_W \otimes \epsilon_W$.  Applying this map to the tensor product of the meaning vectors of the words provides us with the following meaning vector for the sentence:
\begin{eqnarray*}
\overrightarrow{\mbox{men kill cute dogs}} & = & \left(\epsilon_W \otimes 1_S \otimes \epsilon_W \otimes \epsilon_W\right) (\overrightarrow{\text{men}} \otimes \overrightarrow{\text{kill}} \otimes \overrightarrow{\text{cute}} \otimes \overrightarrow{\text{dogs}})\\
&=&   \sum_{ijk} \ c_{ijk}  \ \langle \overrightarrow{men} \mid \overrightarrow{w}_i\rangle  \, \langle   \overrightarrow{w}_k  \mid \overline{\text{cute}}(\overrightarrow{\text{dogs}}) \rangle\,   \overrightarrow{s}_j \\
&=& \sum_{ijk} \sum_{lm}\ c_{ijk} c_{lm} \ \langle \overrightarrow{men} \mid \overrightarrow{w}_i\rangle  \langle\overrightarrow{w}_k \mid \overrightarrow{w}_l\rangle \langle \overrightarrow{w}_m \mid \overrightarrow{dogs}\rangle \ \overrightarrow{s}_j \\
\end{eqnarray*}

As an example of a more complicated sentence, consider `men do not kill dogs'. For the  meaning vectors of the auxiliary  `do',  the infinitive  `kill', and the  negation word   `not', we have the following \cite{CCS,PrellerSadr}:
\[
(\overrightarrow{\text{kill}} \in W \otimes S \otimes W, \sigma^r\cdot j \cdot n^l), \quad 
(\overrightarrow{\text{do}} \in W\otimes S \otimes S \otimes W, n^r \cdot s \cdot j^l \cdot \sigma), \quad
(\overrightarrow{\text{not}} \in W \otimes S \otimes S \otimes  W, \sigma^r \cdot j \cdot j^l \cdot \sigma)
\] 
These are obtained by making the semantic assumption that the vector spaces for  $j$ and $s$ are the same, that is $S$, and that the vector spaces for $n$ and $\sigma$ are the same, that is $W$. The first assumption is justified by the fact that transitive and infinitive transitive verbs both have the same meaning, hence live in the same vector space. The second assumption is justified by the fact that the $\sigma$'s are place holders for the subject.  The pregroup morphism corresponding to the grammatical reduction of the  sentence is as follows:
\[
\alpha \quad = \quad (1_s \cdot \epsilon^l_j \cdot \epsilon^l_{j}) \circ (\epsilon^r_n \cdot 1_s \cdot 1_{j^l} \cdot \epsilon^r_{\sigma} \cdot 1_j \cdot 1_{j^l} \cdot \epsilon^r_{\sigma} \cdot 1_j \cdot \epsilon^l_n)
\]
The linear map $f$ corresponding to the above is:
\[
 (1_S \otimes \epsilon_S \otimes \epsilon_S) \circ (\epsilon_W \otimes 1_S \otimes 1_{S} \otimes \epsilon_W \otimes 1_S \otimes 1_{S} \otimes \epsilon_W \otimes 1_J \otimes \epsilon_W)
\]
The meaning of `not' is generated from a linear map $\overline{\text{not}}: S \to S$, which takes a sentence and negates it. The meaning   of `do' is generated  from the identity $1_S$. These are as follows:
\[
\overrightarrow{\text{do}} \ := \  \left (1_W\otimes \left((1_S \otimes 1_S) \circ \eta_S \right) \otimes 1_W\right) \circ  \eta_W
\hspace{2cm}
\overrightarrow{\text{not}} \ := \  \left (1_W \otimes \left((1_S \otimes \overline{\text{not}}) \circ \eta_S \right) \otimes 1_W\right) \circ  \eta_W
\] 
The corresponding computations provide the following  meaning vector for the  sentence~\cite{CCS,PrellerSadr}: 
\[
\overline{\text{not}} \, \Big(\overrightarrow{\mbox{men kill dogs}}\Big)
\]
This is the \emph{negation} of the meaning of the positive version of the sentence.

So far we have not said anything about the concrete shape of $ W$ and  $S$. Moreover, for the general setting to work, we need concrete vectors for words with compound types whose vectors live in tensor spaces and cannot be directly obtained from distributional models. For example,  the verb lives in $W \otimes S \otimes W$ and has to act on the subject and object, the adjective lives in $W \otimes W$ and has to modify the noun,  the negation word $\overline{\text{not}}$ is a linear map that acts on the abstract   space $S$ and  has to negate the meaning of a sentence and so on.  Previous work presented some solutions for the limit truth-theoretic case,  where $S$ was the two dimensional space of $\mid\!\text{1}\rangle$ and $\mid\!\text{0}\rangle$, and $\overline{\text{not}}$ was taken to be the linear map corresponding to  the matrix of the swap gate~\cite{CCS}. Later,  $S$ was also taken to be a  high dimensional vector spaces concretely built from the distributional models \cite{GS,GS2}. The next section summarises these  results. 

%For instance, a truth-theoretic meaning is obtained by taking   $S$  to be the 2-dimensional space with basis vectors true $\mid\!\!1\rangle$ and false $\mid\!\!0\rangle$. A corresponding definition for the vector of the verb --- that it is $\mid\!\!1\rangle$ if dogs chase cats and is $\mid\!\!0\rangle$ otherwise--- and our  definition for the vector of the sentence ensure that these values propagate to the vector of the sentence. So $\overrightarrow{\text{dogs} \ \text{chase} \ \text{cats}}$ becomes $\mid\!\!1\rangle$ whenever `dogs chase cats' is true, and $\mid\!\!0\rangle$ otherwise. 
%
%

\section{Concrete Spaces and Experimental Results}\label{impl}
In previous work,  we provided  an algorithm  for  building vectors for words with compound types and instantiated it  for the particular case of  intransitive and transitive sentences with simple subjects and objects~\cite{GS,GS2}. Work in progress extends these cases to adjective-noun phrases in transitive sentences   as subject and object. We  evaluated the resulting vectors on a disambiguation task performed  on real data obtained from the British National Corpus (BNC). This corpus  consists of about 6 million sentences and 500,000 unique lemmas.

\subsection{Concrete Constructions}
We build a vector space $N$ from the BNC  by taking its 2000 most occurring lemmas as basis vectors $\overrightarrow{n}_i$; this restriction is both for computational purposes and also to be able to compare our results to related work~\cite{Lapata}.  Vectors of  $N$ are built by counting co-occurence and normalizing by TF/IDF. We take $W$ to be $N$ and $S$ to be $N \oplus (N \otimes N) \oplus (N \otimes N \otimes N)$. In the latter, the first $N$ encodes meanings of intransitive sentences; these are unary relations; $N \otimes N$ encodes meanings of transitive sentences, which are binary relations, and $N \otimes N \otimes N$ is for meanings of ditransitive sentences, which are ternary relations. These relations are denoted by weighted sets of singletons, pairs, and triples, respectively. The set of singletons for an intransitive verb denotes the subjects to which the verb has been applied throughout the corpus. The set of pairs denotes the subject-object pairs that have been related by a transitive verb, etc. The weights represent  \emph{the extent according to which} the verb has acted on or related the nouns; these  are learnt from the corpus  in the following two ways. These are instantiations of  a general procedure for building vectors for a word of any compound type~\cite{GS}. 

\begin{enumerate}
\item {\bf Categorical (1)}. The  meaning vector of an intransitive verb $\overrightarrow{v} \in N$ is $\sum_i c_i \overrightarrow{n}_i$, where each $c_i$ is obtained by summing the vectors of all the subjects of $\overrightarrow{v}$ throughout the BNC. The  meaning vector of a  transitive verb $\overrightarrow{tv} \in N \otimes N$ is $\sum_{ij} c_{ij} \overrightarrow{n}_i \otimes \overrightarrow{n}_j$, where each $c_{ij}$ is the sum of the tensor products of all the subjects and objects that $tv$ has related in the BNC, and similarly for  ditransitive verbs.  

\item {\bf Categorical (2)}.  Here we simply take the Kronecker products of the context vectors of the verbs, abusing the notation we still refer to it as $\overrightarrow{v}$. So for the intransitive verb, this is the vector obtained by a normalized count of co-occurence $\overrightarrow{v} \in N$. For the transitive verb it is $\overrightarrow{v} \otimes \overrightarrow{v} \in N\otimes N$ and for the ditransitive it is $\overrightarrow{v} \otimes \overrightarrow{v} \otimes \overrightarrow{v} \in N \otimes N \otimes N$. 
\end{enumerate}

The above vectors are embedded in the spaces prescribed for them by a DisCoCat, that is $N \otimes S$ for $\overrightarrow{v}$ and $N \otimes S \otimes N$ for $\overrightarrow{tv}$,  by diagonalisation, that is by padding the  non-diagonal elements by zeros.

%For instance, suppose $N$ has only two bases, hence the relational vector of the verb will be a $2 \times 2$ vector and its typological vector is a $4 \times 4$ one.  Suppose further that the transitive verb `kill' has only occurred twice in the corpus, in sentences  `s$_1$ kill o$_1$' and `s$_2$ kill o$_2$'. Then the relational vector of the kill will be $\text{s}_1 \otimes \text{s}_2 + \text{s}_2 \otimes \text{o}_2$ and lives in $N \otimes N$. The typological vector will be a diagonal $4 \times 4$ matrix with the relational vector as its diagonal.  

\subsection{Evaluation Task, Dataset, and Results}
We evaluated these concrete methods on a disambiguation task. The general idea behind this  task is that some verbs have more than one meaning  and the sentences in which they appear  disambiguate them. Suppose a verb $v$  has two meanings $a$ and $b$ and we want to decide whether $v$ means $a$ or $b$  whenever it occurs in a sentence $s$.  To implement this task, we  chose 10 ambiguous  transitive verbs from the most frequent verbs of the BNC. For each verb, two different non-overlapping meanings were retrieved via the JCN measure of  information content synonymy applied to  WordNet  synsets \cite{jiang}. For instance  one of our chosen verbs was `meet', for which we obtained meaning $a$:`visit' and meaning $b$:`satisfy'. For each original verb, ten sentences containing that verb  (in the same role) were retrieved from the BNC; for example,  one such  sentence  for $v$: `meet' is  $s$: `the system met the criterion'.  For each such sentence, we generated two other  sentences by substituting the verbs of those sentences by $a$ and $b$, respectively.  For instance,  `the system satisfied the criterion' and `the system visited the criterion' were generated for the first meaning of `meet'. This procedure provided us with 200 pairs of sentences. Note that the generated sentences only make sense for the correct meaning of the verb; for instance, `the systems satisfied the criterion' does make sense, whereas `the system visited the criterion' does not.  The goal is to verify that  the sentences `the system met the criterion' and `the system satisfied the criterion'  have a high degree of semantic similarity, whereas the  sentences  `the system met the criterion' and `the system visited  the criterion' have a  low degree of similarity. The result of this verification  disambiguates the verb. For the case of  `meet',  we are verifying that  it  means `satisfy' (and not `visit') in the sentence `the system met the criterion'.

We experimented with two datasets, one for intransitive sentences from~\cite{Lapata} and an extension of it for transitive sentences, as described above (the extension to sentences with adjective0noun phrases and subject and object is straightforward and preserves the results). An example of the second dataset is provided  in Table 3. 

\begin{table}[h]
\begin{center}
 \begin{tabular}{|c|c|c|}
 \hline
 &Sentence 1 & Sentence 2\\
 \hline
 \hline
 
1& the system met the criterion & the system visited the criterion  \\
\hline
2& the system met the criterion  & the system satisfied the criterion\\
\hline
3& the child met the  house & the child visited the  house\\
 \hline
 4& the child met the  house  & the child satisfied the  house \\
 \hline
5&  the child showed  interest  & the child pictured  interest \\
 \hline
 6&  the child showed  interest  & the child expressed  interest \\
\hline
7& the map showed the location & the map pictured the location\\
\hline
8& the map showed the location & the map expressed the location\\ 
\hline
\end{tabular}
\end{center}
\label{sampledataset3}
\caption{Sample Sentence Pairs from the Transitive  Dataset.}
\end{table}

We built vectors for  nouns by using the usual distributional co-occurence method. Then built vectors for verbs using both of the methods described above, and finally built vectors for sentences using the DisCoCat prescription. After the sentence vectors were built, we measured the similarity of each  pair using the cosine of their angles  on the scale of the real numbers in $[0,1]$. In order to  judge the performance of our method, we followed guidelines from related work~\cite{Lapata}.  We distributed our data set among  25 volunteers who were asked to  rank each pair based on how similar they thought they were. The ranking was between 1 and 7, where 1 was almost dissimilar and 7 almost identical. To be in line with related work~\cite{Lapata}, each pair was also given a HIGH or LOW classification by us.  However, these scores are solely based on   our personal judgements and on their own they do not provide a very reliable measure of comparison.  The correlation of the model's similarity judgements with the human judgements was calculated using Spearman's $\rho$. This is  a rank correlation coefficient ranging from -1 to 1 and  provides a  more robust metric (in contrast with the LOW/HIGH metric)  by which models are  ranked and compared.  

%It is assumed that inter-annotator agreement provides the theoretical maximum $\rho$ for any model for this experiment, hence these are used as the \emph{UpperBound}. Taking the cosine measure of the verb vectors while ignoring the noun was taken as the \emph{Baseline}.

\begin{table}[h]
\begin{center}
\begin{tabular}{|lll|c|}
\hline
Model & High & \quad Low & \quad $\rho$\\
\hline
\hline
Baseline & 0.47 & \quad 0.44 & \quad 0.16\\
\hline
\hline
Add & 0.90 & \quad 0.90 & \quad 0.05\\
Multiply & 0.67 & \quad 0.59 & \quad 0.17 \\
\textbf{Categorical (1)} & \textbf{0.73} & \quad\textbf{0.72} & \quad \textbf{0.21}\\
\textbf{Categorical (2)} & \textbf{0.34} & \quad\textbf{0.26} & \quad \textbf{0.28}\\
\hline
\hline
UpperBound & 4.80 & \quad 2.49 & \quad 0.62 \\
\hline
\end{tabular}
\end{center}
\label{results4}
\caption{Results of the 1st and 2nd Compositional Disambiguation Experiments.}
\end{table}

The results  of these calculations for our datasets are presented in Table 4.  
The additive and multiplicative rows have, as composition operation, vector addition and component-wise multiplication.
The \emph{Baseline} is from a non-compositional approach; it is obtained by comparing the verb vectors of each pair directly and ignoring their subjects and objects. The \emph{UpperBound} is set to be inter-annotator agreement. According to the  $\rho$ measure,   both of  our methods outperform the other methods.

%The increase of $\rho$ reflects the compositionality  of our model:  its performance increases with the increase in syntactic complexity.  Based on this, we would like to believe  that more complex datasets and experiments which for example include adjectives and adverbs  shall lead to even better results.  

\section{A Compositional Distributional Model of Meaning on  Lambek Monoids} 

%As it is the case in categorical logic, $\leq$ of a partial order algebra becomes a $hom$ object  in a category.  The categorical analogues to all of the above partial order algebras, each referred to as a type-logic category ${\cal T}$, have the appropriate structural morphisms. For instance, each residuated monoid is a closed monoidal category. Over such a category, we define a functor that interprets it in the category of finite dimensional vector spaces. 

The contribution of this paper is to extend the DisCoCat model from pregroups to  Lambek monoids. 
Rather than pairing $FVect$ with a  monoidal category, we will rely on  a functorial passage from a monoidal category  to $FVect$.  There are two reasons  to opt for a  functorial passage. Firstly, it makes a closer connection to the original  semantic models of natural language \cite{montaguebook}, which were based on a  homomorphic passage from sentence formation rules  to set theoretic operations. In our model, this homomorphic  passage is formalised via a functor between  categories of grammatical  reductions and meaning. Contrary to those models  \cite{montaguebook}, however, our model is not limited to sets and set-theoretic operations and is generalised  to vectors and vector composition operations. This brings us to our second reason: a homomorphic passage to the category of vector spaces is not a one-off development especially tailored for our purposes. It is   an example of a more general construction, namely,  a passage long-known in Topological Quantum Field Theory (TQFT). This general passage was first developed in   \cite{Atiyah} in the context of TQFT and was given the name `quantisation', %as it provided a quantitative counterpart (i.e.~in terms of vectors) to qualitative structures that expressed forms (in terms of cobordisms and manifolds).  
 as it adjoins `quantum structure' (in terms of vectors) to a purely topological entity, namely the cobordisms representing the topology of manifolds.
Later, this passage was generalised to abstract mathematical structures and recast in terms of functors whose co-domain was $FVect$~ \cite{BaezDolan, Kock}.  This is exactly what is happening in our semantic framework: the sentence formation rules are formalised using type-logics and assigned quantitative values in terms of vector composition operations. This procedure makes our passage from grammatical structure to vector space meaning  a  `quantisation' functor. Hence, one can say that what we are developing  here is a grammatical quantum field theory for Lambek monoids. The detailed constructions of this paper can be worked out in a similar fashion for the  pregroup-based framework of previous work \cite{CCS,CSC}.

\subsection{Monoidal Bi-Closed Categories}
A Lambek monoid is a monoidal bi-closed category. Similar to the case of pregroups,  its objects are the elements of the partial order and its morphisms are  provided by the ordering relation. The  monoid multiplication is a non-symmetric tensor, whose residuals are the left and right implications. To see this, we recall some definitions.

%\begin{definition}[closed monoidal category]
%...
%\end{definition}

%\begin{definition}[compact (closed) category]
%...
%\end{definition}

%The following are then straightforward.

%\begin{proposition}
%Each compact (closed) category is a closed monoidal category.
%\end{proposition}

%\begin{proposition}
%Each pregroup is a compact (closed) category.
%\end{proposition}

%\begin{proposition}
%Each residuated monoid is a closed monoidal category.
%\end{proposition}

%\paragraph{Graphical language for symmetric monoidal categories \cite{Joyal-Street}.}

%\paragraph{Graphical language for closed monoidal categories \cite{Baez}.}

%\paragraph{Graphical language for compact (closed) categories \cite{Kelly}.}

%\smallskip
%\begin{definition}
A monoidal bi-closed category is a category with a monoidal tensor $\otimes$ and its unit $I$, such that for all pairs of objects $A,B$ of the category,   there exist a pair of  objects $A \multimap B$, $A \multimapinv B$ and a pair of morphisms  as follows: 
\[
ev^l_{A,B} \colon A \otimes (A\multimap B) \to B \hspace{2cm}
ev^r_{A,B} \colon (A \multimapinv B) \otimes B \to A
\]
These morphisms  are referred to  as \emph{left} and \emph{right evaluations}. They are  such that for every pair of arrows $f\colon (A \otimes C) \to B$ and $g \colon (C \otimes B) \to A$  
there  exist two unique  morphisms  as follows:
\[
\Lambda^l(f) \colon C \to A \multimap B \hspace{2cm} \Lambda^r(g) \colon C \to A \multimapinv B
\]
These morphisms are referred to as \emph{left and right currying}; they make the following diagrams commute:
\begin{diagram}
A \otimes C &\rTo^{1_A \otimes \Lambda^l(f)} &A \otimes (A \multimap B)
\hspace{2cm} & C\otimes B &\rTo^{\Lambda^r(g)\otimes 1_B}&(A\multimapinv B) \otimes B\\
&\rdTo_{f}&\dTo_{ev^l_{A,B}} 
\hspace{2cm} &&\rdTo_{g} &\dTo_{ev^r_{A,B}}\\
&&\qquad \qquad \quad B\hspace{2cm}& &&A
\end{diagram}
%\end{definition}

\medskip
Given  a morphism $f\colon A \to B$, its \emph{name} $\ulcorner f \urcorner \colon I \rTo^{\Lambda^l(f)} A \multimap B$, is obtained by  currying the morphism $(A\otimes I) \rTo^{f} B$ (dropping the precomposition with the  $\iota_A$ morphism, i.e. $A \otimes I \rTo^{\iota_A} A$). In  order to be coherent with the above left-right notation, we define a left and a right name for $f$, obtained by currying the morphisms  $(A \otimes I) \rTo^{f} B$  and $(I \otimes A) \rTo^{f} B$,   as follows:
\[
\ulcorner f \urcorner^l \colon I \rTo^{\Lambda^l(f)} A \multimap B
\hspace{2cm}
\ulcorner f \urcorner^r \colon  I \rTo^{\Lambda^r(f)} B \lto A
\]
Evaluating these names makes the following diagrams commute:
\begin{diagram}
A \otimes I &\rTo^{1_A \otimes \ulcorner f \urcorner^l}  & A \otimes (A\multimap B) \hspace{2cm} & 
I \otimes A &\rTo^{\ulcorner f \urcorner^r \otimes 1_A}  & (B \lto A) \otimes A
\\
&\rdTo_{f}&\dTo_{ev^l_{A,B}} 
\hspace{2cm} &&\rdTo_{f} &\dTo_{ev^r_{A,B}}\\
&&\qquad \qquad \quad  B\hspace{2cm}& &&B
\end{diagram}
In other words, we have the following two equations:
\[
ev^l_{A,B} \circ (1_A \otimes\ulcorner f \urcorner^l) \quad =  \quad f
\hspace{2cm}
ev^r_{A,B} \circ (\ulcorner f \urcorner^r  \otimes 1_A) \quad =  \quad f
\]
The names of the  identity  morphism $1_A \colon A \to A$ are obtained as special cases of the above construction. These are defined as follows:
\[
\ulcorner 1_A \urcorner^l \colon I \rTo^{\Lambda^l(1_A)} A \multimap A 
\hspace{2cm}
\ulcorner 1_A \urcorner^r \colon  I \rTo^{\Lambda^r(1_A)} A \lto A
\]
Evaluating the above  leads to two similar commutative diagrams; the  equations corresponding to these  diagrams  provide us with  the monoidal version of \emph{yanking}:
\[
ev^l_{A,A} \circ (1_A \otimes\ulcorner 1_A \urcorner^l) \quad = \quad 1_A
\hspace{2cm}
ev^r_{A,A} \circ (\ulcorner 1_A \urcorner^r  \otimes 1_A) \quad = \quad 1_A
\]

\subsection{A Quantisation Functor for Lambek Monoids}\label{quantfunct}
The quantisation functor  preserves all the syntactic structure but ``forgets'' the order of the words. This order will be taken care of in our explicit definition of the meaning vector of a sentence in Definition~\ref{disco}.  

\medskip
\begin{definition}\label{quantdef}
Given a Lambek monoid  ${\cal L}$ and the category of finite dimensional vector spaces $FVect$ over  $\mathbb{R}$, the quantisation functor ${\cal Q}\colon {\cal L} \to FVect$ is a strongly monoidal functor,  satisfying the following:  
\begin{eqnarray}
{\cal Q}(1) &:= & \ \mathbb{R}\\
{\cal Q}(a \cdot b) \ &\cong&   \ {\cal Q}(b \cdot a) \quad \ \ \ :=  \quad {\cal Q}(a) \otimes {\cal Q}(b)\\
{\cal Q}(a \multimap b) \ &\cong& \ {\cal Q}(b \multimapinv a)  \quad := \quad{\cal Q}(a) \Rightarrow {\cal Q}(b) 
%{\cal Q}(\pi \multimapinv \pi') & \cong& {\cal Q}(\pi') \Rightarrow {\cal Q}(\pi)
\end{eqnarray}
\end{definition} 
The tensor product  in $FVect$ is symmetric, so  (4)  forgets the order and (5) collapses the two implications of $\cal L$  to just one. Also, $FVect$  only has one evaluation  and  currying map as follows: 
\[
ev_{A,B \colon }A \otimes (A \Rightarrow B) \to B \hspace{2cm}
\Lambda(f) \colon C \to (A \Rightarrow B), \quad \text{for} \ f\colon (A \otimes C) \to B
\]

The closed object $V \Rightarrow W$ is the set of all linear maps from $V$ to $W$, made into a vector space in the standard way, via the isomorphism  $V \Rightarrow W \cong V^* \otimes W
%, see~\cite{Baez}.  We use vector spaces that are on the field of reals and have a fixed basis, hence the isomorphism $V^* \cong V$ becomes canonical and we obtain  $V \Rightarrow W \cong 
 \cong V \otimes W$.  Hence,  the closed structure of $FVect$ determines its compact structure. Given this fact,  the evaluation map of $FVect$ is the same as the corresponding co-unit of the compact adjunction. Hence,  for $V,W$ respectively spanned by $\{\overrightarrow{v}_i\}_i, \{\overrightarrow{w}_j\}_j$, and a vector  $\overrightarrow{v} \in V$, the morphism 
\[
V \otimes (V \Rightarrow W) \rTo^{ev_{V,W}} W
\]
is as in the compact closed case, that is:
\[
ev_{V,W} \left(\overrightarrow{v}\otimes \left(\sum_{ij} c_{ij}\  \overrightarrow{v}_i \otimes \overrightarrow{w}_j\right) \right) \quad =  \quad
\sum_{ij} c_{ij}\  \langle \overrightarrow{v} \mid \overrightarrow{v}_i\rangle \ \overrightarrow{w}_j 
\]

The quantisation of  atomic types yields atomic vector spaces:    for a noun phrase $n$ we  stipulate   ${\cal Q}(n) :=  N$  and  for a declarative sentence $s$  we stipulate ${\cal Q}(s) :=  S$. The quantisation  of compound types yields closed  vector spaces.  The quantisation of  an intransitive verb, which has the  type $n \multimap s$, is computed as follows:
\[
{\cal Q}(n \multimap s) \quad \cong \quad {\cal Q}(n) \Rightarrow {\cal Q}(s) \quad := \quad N \Rightarrow S
\]
The quantisation of an adjective, which has the type $n \multimapinv n$,  is computed as follows:
\[
{\cal Q}(n \multimapinv n)  \quad \cong \quad {\cal Q}(n) \Rightarrow {\cal Q}(n)\quad := \quad N \Rightarrow N
\]
The quantisation of a transitive verb, which has the type $(n \multimap s) \multimapinv n$,   is computed as follows:
\[
{\cal Q}((n \multimap s) \multimapinv n) \quad \cong \quad {\cal Q}(n) \Rightarrow ({\cal Q}(n) \Rightarrow {\cal Q}(s)) \quad := \quad  N \Rightarrow (N \Rightarrow S)
\]

%\noindent
%This definition is only partly categorical. It can be re-notated to something that is more categorical, by replacing the vectors with their initial map, e.g.~$\overrightarrow{v} \in V$ with $\mathbb{R} \rTo^{\overrightarrow{v}} V$. So the above definition will look like
%\[
%I \rTo^{\ \overrightarrow{\text{w}_1 \text{w}_2 \cdots \text{w}_n}\ } S \quad := \quad
%{\cal Q}(\alpha) \left((I\rTo^{\ \overrightarrow{\text{w}_1}\ }W_1)  \otimes (I \rTo^{\ \overrightarrow{\text{w}_2}\ } W_2) \otimes \cdots \otimes (I \rTo^{\ \overrightarrow{\text{w}_n}\ } W_n)\right)
%\]
%I am not sure what to think of equality of morphisms, may be the above is read as: a  morphism  $I \rTo S$ is obtained by tensoring a set of morphisms $I \rTo W_i$,  then applying to them ${\cal Q}(\alpha)$.

Monoidal  meaning vectors of strings of words are computed by applying  the quantisation  of their grammatical reduction map to the tensor products of the quantisations of their words. More precisely, we have:

\smallskip
\begin{definition}\label{disco}
The monoidal  meaning vector of a string  $\text{w}_1 \text{w}_2 \cdots \text{w}_n$ consisting of $n$ words is: 
\[
\overrightarrow{\text{w}_1 \text{w}_2 \cdots \text{w}_n} \quad := \quad
{\cal Q}(f) \left(\overrightarrow{\text{w}_1} \otimes \overrightarrow{\text{w}_2} \otimes \cdots \otimes  \overrightarrow{\text{w}_n}\right)
\]
where for $t_i$ the grammatical type of the word $\text{w}_i$, the map $f$ is the monoidal grammatical reduction map of the string, that is   $t_1 \cdot t_2 \cdot \  \cdots\  \cdot t_n \rTo^{f} s$    and ${\cal Q}(f)$ is defined as follows:
\[
{\cal Q}(t_1 \cdot t_2 \cdot \  \cdots\  \cdot t_n \rTo^{f} s) \quad := \quad
{\cal Q}(t_1 \cdot t_2 \cdot \cdots \cdot t_n) \  \rTo^{{\cal Q}(f)}  \ {\cal Q}(s) 
\] 
\end{definition}

For example, the  ${\cal Q}(f)$ of  an intransitive sentence with $f  = n \cdot (n \multimap s) \rTo^{ev^l_{n,s}} s$, is   computed as follows:
\begin{eqnarray*}
{\cal Q}\left(n \cdot (n \multimap s) \rTo^{ev^l_{n,s}} s\right) &\quad \cong \quad&
{\cal Q}\left(n \cdot (n \multimap s)\right) \rTo^{{\cal Q}\left(ev^l_{n,s}\right)} Q(s)\\
&\cong& N \otimes (N \Rightarrow S) \rTo^{ev_{N,S}} S
\end{eqnarray*}
The monoidal meaning of an intransitive sentence `men kill' is as follows
\[
\overrightarrow{\text{men \ kill}}  := {\cal Q}(f)\left(\overrightarrow{\text{men}} \otimes \overrightarrow{\text{kill}}\right) \\
\]
Given that a vector $\overrightarrow{v}$ in a vector space $V$ is represented by the morphism $I \rTo^{\overrightarrow{v}} V$, the  details of this computation  become as follows:
\begin{eqnarray*}
\overrightarrow{\text{men \ kill}}  &\cong &  ev_{N,S}  \left(\left(I \rTo^{\overrightarrow{\text{men}}} N\right) \ \otimes \ \left (I \rTo^{\overrightarrow{\text{kill}}} (N \Rightarrow S)\right) \right)\\
&=&  ev_{N,S}  \left(I \otimes I \rTo^{\overrightarrow{\text{men}} \otimes \overrightarrow{\text{kill}}} N \otimes (N \Rightarrow S)\right)\\
&\cong& I \rTo^{\overrightarrow{\text{men}} \otimes \overrightarrow{\text{kill}}} N \otimes (N \Rightarrow S) \rTo^{ev_{N,S}} S \\
&=& I \rTo^{ev_{N,S} (\overrightarrow{\text{men}} \otimes \overrightarrow{\text{kill}})} S
\end{eqnarray*}
This morphism picks a sentence vector from $S$;  given that   the monoidal and compact meanings of the  verb `kill' coincide and both are representable by $\sum_{ij} c_{ij}\  \overrightarrow{n}_i \otimes \overrightarrow{s_j}$,  the above morphism picks  the vector  $\sum_{ij} c_{ij}\  \langle \overrightarrow{\text{men}} \mid \overrightarrow{n}_i\rangle \ \overrightarrow{s_j}$, which is the same vector as the compact meaning vector of `men kill'.

For the transitive  sentence,  the grammatical reduction map is the composition of two  maps as follows:
\[
n \cdot ((n \multimap s)\multimapinv n) \cdot n \rTo^{1_{n} \cdot ev^r_{n, n \multimap s}} n \cdot (n \multimap s) \rTo^{ev^l_{n,s}} s
\]
Hence, $f$ is  the composition $ev^l_{n,s} \ \circ \ \big(1_n \cdot ev^r_{n, n\multimap s}\big)$. The quantisation of this composition   is the composition of quantisations, computed as follows:
\begin{align*}
{\cal Q}\Big(ev^l_{n,s} \ \circ \ \big(1_n \cdot ev^r_{n, n\multimap s}\big)\Big) &\cong {\cal Q}\Big(ev^l_{n,s}\Big) \ \circ \ {\cal Q}\Big(1_n \cdot ev^r_{n, n\multimap s}\Big)  \hspace{2cm} \mbox{by functoriality of} \ {\cal Q}\\
&\cong ev_{N,S}  \ \circ \ \Big({\cal Q}\big(1_n\big) \otimes {\cal Q}\big(ev^r_{n, n\multimap s}\big)\Big) \hspace{1.5cm} \mbox{by monoidality of}\ {\cal Q}\\
&\cong ev_{N,S}  \ \circ \ \Big(1_{{\cal Q}(n)} \otimes ev_{N, N \Rightarrow S}\Big) \hspace{2cm} \mbox{by functoriality of} \ {\cal Q}\\
&\cong ev_{N,S} \ \circ \ \Big(1_N \otimes ev_{N, N\Rightarrow S}\Big) 
\end{align*}
And based on the above, ${\cal Q}(f)$  is computed as follows   for  a transitive sentence:
\begin{align*}
{\cal Q}\left(n \cdot ((n \multimap s)\multimapinv n) \cdot n \rTo^{ev^l_{n,s} \circ \big(1_n \cdot ev^r_{n, n\multimap s}\big)} s\right)& \quad \cong \quad
N \otimes (N \Rightarrow (N \Rightarrow S)) \otimes N \rTo^{ev_{N,S} \circ \big(1_N \otimes ev_{N, N\Rightarrow S}\big)} S
\end{align*}
The monoidal   meaning  vector of a transitive sentence such as  `men kill dogs' simplifies to the following: 
\[
\overrightarrow{\mbox{men kill dogs}} \quad := \quad \left(ev_{N,S} \circ \big(1_N \otimes ev_{N, N\Rightarrow S}\big)\right) 
\left(\overrightarrow{men} \otimes \overrightarrow{kill} \otimes \overrightarrow{dogs}\right)
\]
Following a computation similar to that of the intransitive sentence, the above picks the sentence vector $\sum_{ijk} c_{ijk} \ \langle \overrightarrow{\text{men}} \mid \overrightarrow{n}_i\rangle  
\langle \overrightarrow{n}_k\mid \overrightarrow{\text{dogs}}\rangle\ \overrightarrow{n}_j $ from $S$; this vector is the same  vector as the one obtained in the compact closed case.

For an adjective-noun phrase, $f$ is $(n \multimapinv n) \cdot n \rTo^{ev^r_{n,n}} n$, hence ${\cal Q}(f)$ becomes as follows:
\[
{\cal Q}\left((n \multimapinv n) \cdot n \rTo^{ev^r_{n,n}} n\right) \quad \cong \quad 
 (N \Rightarrow N) \otimes N \rTo^{ \ ev_{N,N}} N
\]
%The monoidal meaning vector of an adjective-noun phrase is computed as follows:
%\begin{eqnarray*}
%\overrightarrow{\mbox{adj noun}} & := & ev_{N,N}\left(\overrightarrow{\text{adj}} \, \otimes \, \overrightarrow{\text{noun}}\right)\\
%&\cong& ev_{N,N} \left((N \Rightarrow N) \otimes N\right)\\
%&\cong& N
%\end{eqnarray*}
whose concrete vector  is the same as in the compact closed case, that is $\sum_{lm} c_{lm}\  \langle \overrightarrow{\text{noun}} \mid \overrightarrow{n}_l\rangle \ \overrightarrow{n}_m$.

For a sentence  with an adjective-noun phrase, such as `men kill cute dogs', ${\cal Q}(f)$ is computed as follows:
\begin{align*}
&{\cal Q}\left(n \cdot ((n \multimap s)\multimapinv n) \cdot (n \lto n) \cdot n 
\rTo^{ev^l_{n,s} \circ \big(1_n \cdot ev^r_{n, n\multimap s}\big) \circ \big(1_{n \cdot (n\multimap n) \lto s} \cdot ev^r_{n,n}\big)} 
s\right)\cong \\
&N \otimes (N \Rightarrow (N \Rightarrow S)) \otimes (N \Rightarrow N) \otimes N 
\rTo^{ev_{N,S} \circ \big(1_N \otimes ev_{N, N\Rightarrow S}\big) \circ \big(1_{N \otimes (N \Rightarrow (N \Rightarrow S))} \otimes ev_{N,N}\big)} 
S
\end{align*}
The computation of the monoidal meaning vector of the above is done by substituting the meaning vector of the adjective-noun phrase for the meaning of the object in the meaning vector of the transitive sentence. It provides us with the same meaning vector as in the compact closed case.

To compute the monoidal vector meaning of the sentence `men do not kill cute dogs', we take  `do' to be the name of the identity  morphism on $N\Rightarrow S$ and  `not' to be  the name of $\overline{\text{not}}$  on the morphism $(N\Rightarrow S) \to (N \Rightarrow S)$. That is:
 \begin{eqnarray*}
 \overrightarrow{\text{do}} &:=&  \ulcorner 1_{N \Rightarrow S}\urcorner \colon \mathbb{R} \to ((N\Rightarrow S) \Rightarrow (N\Rightarrow S))\\
    \overrightarrow{\text{not}} &:=&  \ulcorner \overline{\text{not}}\urcorner \colon \mathbb{R} \to ((N\Rightarrow S) \Rightarrow (N\Rightarrow S))
    \end{eqnarray*}
The meaning vector of the sentence will  then simplify to the following:
\[
 ev_{N,S}  \left(\overrightarrow{men} \otimes \overline{\text{not}}(\overrightarrow{kill}(-, \overrightarrow{dogs}))\right)
\]
So the monoidal meaning of a negative sentence is obtained by first applying the meaning of the verb to the object,  then negating it, then applying the result to  the subject.  This procedure  is different from the compact closed case, which was obtained by first applying the meaning of the verb to the meanings of subject and object, then negating it.  The monoidal meaning results in the following vector:
\[
\left\langle \overrightarrow{\text{men}} \left | \right. \overline{\text{not}} \left(\sum_{ijk} c_{ijk}\  \overrightarrow{n}_i  \langle \overrightarrow{n}_k \mid\overrightarrow{\text{dogs}}\rangle \ \overrightarrow{s}_j\right)\right \rangle
\]
 Whereas the compact meaning results in the following vector: 
 \[
  \overline{\text{not}} \left(\sum_{ijk} c_{ijk}\ \langle \overrightarrow{\text{men}} \mid \overrightarrow{n}_i \rangle \langle \overrightarrow{n}_k \mid\overrightarrow{\text{dogs}}\rangle \ \overrightarrow{s}_j\right)
 \]
These two vectors  only become equivalent in  special cases, for instance one in  which the   meaning of $\overline{\text{not}}$ can only act on $S$ (and not on $N$). This happens in the   truth theoretic case, in which $N$ has many dimensions,  $S$ has only two  dimensions, and  $\overline{\text{not}}$  swaps the basis vectors, that is, it is the linear map corresponding the the matrix $\left(\begin{array}{cc} 1&0\\0&1\end{array}\right)$. In this case,   $\overline{\text{not}}$ can be applied to any vector in $S$, but it is not defined for the vectors in $N$.  If  $\overline{\text{not}}$ is also defined for $N$, for instance it is a certain  permutation of basis vectors, then in the computation of the monoidal meaning vector of the sentence, it can either apply to the basis elements of $N$, that is to $\overrightarrow{n}_i$, or to the basis elements of $S$, that is to $\overrightarrow{s}_j$. Only the latter will provide a result which is the same as the meaning vector of the sentence in the compact case. 

%\footnote{We suppose that the types have been disambiguated,  using for example the role of the word in the given sentence. This  helps us avoid having to deal with multiple type assignments. For instance the word `fish' might be a noun phrase or a verb.}

% These are specified below.
%\begin{align*}
%\overrightarrow{\text{cute}} \quad&\simeq \quad \ulcorner f \urcorner \colon I \to (N \Rightarrow N)\\
%\overrightarrow{\text{does}} \quad&\simeq \quad \ulcorner 1_{N \Rightarrow S}\urcorner \colon I \to ((N\Rightarrow S) \Rightarrow (N\Rightarrow S))\\
%\overrightarrow{\text{not}} \quad&\simeq \quad \ulcorner g\urcorner \colon I \to ((N\Rightarrow S) \Rightarrow (N\Rightarrow S))
%\end{align*}

\section{Diagrammatic Reasoning in Monoidal Bi-Closed Categories}\label{sec:diagrams}

\usetikzlibrary{arrows,decorations,backgrounds,positioning,fit}

	\tikzset{func/.style={shape=rectangle,rounded corners=8,minimum width=2cm,minimum height=.5cm,draw}}
	\tikzset{claspnode/.style={shape=circle,minimum width=0.25cm,fill=white,draw}}

\begin{figure}[h]
\begin{center}
\epsfig{figure=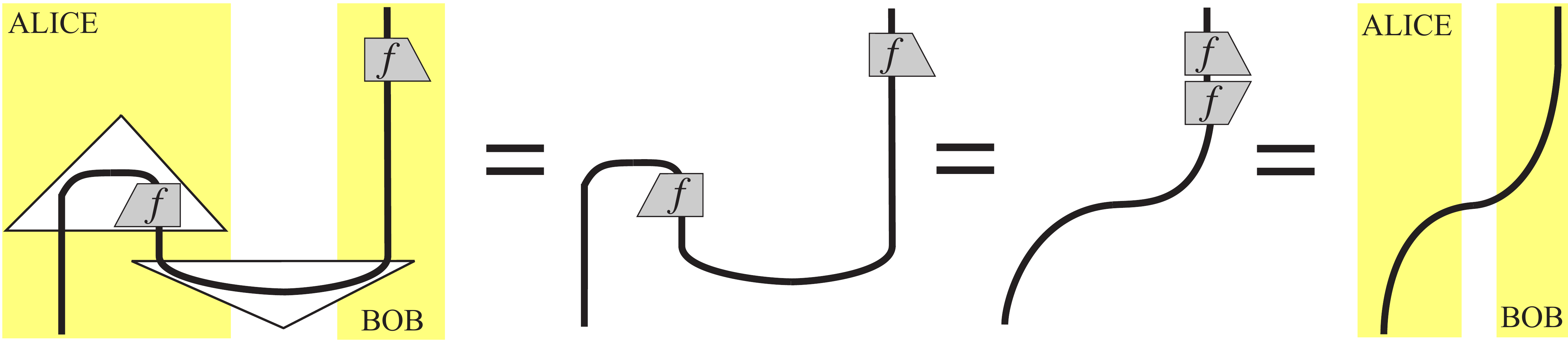,width=280pt}
\end{center}
\caption{Diagram of Information Flow in the Teleportation Protocol}
\label{fig:teleport}
\end{figure}

One of the principal advantages of a DisCoCat is that it is equipped with the sound and complete diagrammatic calculus of compact closed  categories~\cite{Kelly}. We refer the reader for further details to previous work, but in a nutshell, this diagrammatic calculus depicts the information flows that happen among the words of a sentence and which produce a meaning for the full sentence. This information flow happens when certain  objects cancel out via \emph{yanking}; this procedure is depicted by pulling sequences of  connected \emph{cup} and \emph{cap} structures, representing $\epsilon$ and $\eta$ maps,  and turning them into straight lines. If the information flow happens in stages, the process is depicted via equivalences between the resulting diagram of each stage. The equivalences   translate into equations  between categorical morphisms corresponding to each diagram and   provide a symbolic proof of the claim that for instance, a sentence is grammatical. The yanking operations have also been useful to elegantly describe the flow of information in quantum protocols such as teleportation, as depicted in Figure~\ref{fig:teleport}, from~\cite{AbrCoe2}.

\begin{figure}[h!]
\begin{center}
\raisebox{-3mm}{\epsfig{figure=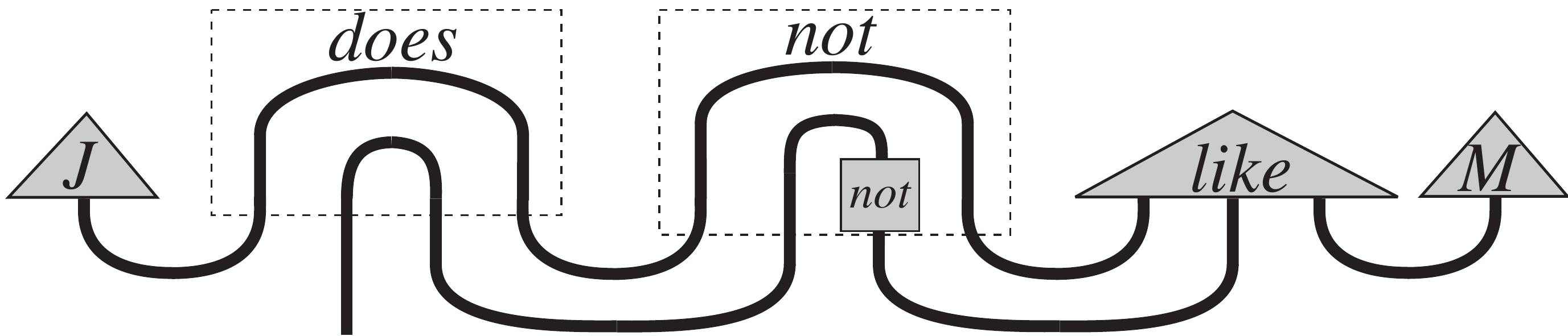,width=280pt}}\vspace{-1mm}
\end{center}
\caption{Diagram of Information Flow in the Negative Transitive Sentence}
\label{fig:neg}
\end{figure}

In a linguistic context, the yanking operations  depict the flow of  information  between words, as depicted in Figure~\ref{fig:neg},  from~\cite{CCS,PrellerSadr}. 

The ability to compute with diagrams  is not lost in the move from compact to %monoidal  and in particular 
monoidal bi-closed categories. A diagrammatic calculus for these categories has been developed by Baez and Stay~\cite{Baez}. Similar to the graphical calculus of compact closed categories, the graphical calculus of monoidal bi-closed categories   can  be applied to both depict the grammatical reduction of Lambek monoids and the flow of information between meanings of the words in a sentence.

\begin{figure}[htbp]
		\centering
			\begin{tikzpicture}[thick,xscale=.75, yscale=0.65]
\hspace{-1cm}
			\tikzset{func/.style={shape=rectangle,rounded corners=8,minimum width=1cm,draw}}
			\tikzset{claspnode/.style={shape=circle,minimum width=0.25cm,fill=white,draw}}

			\draw (-14,2.5) node {$id_A$};
			\draw[->] (-14,2) -- node [left] {$A$} (-14,-1);
\hspace{2cm}
			\draw (-11,2.5) node {$f : A \to B$};
			\node[func] (f) at (-11,.5) {$f$};
			\draw[->] (-11,2) -- node [left] {$A$} (f);
			\draw[->] (f) -- node [left] {$B$} (-11,-1);

			\end{tikzpicture}
		\caption{Basic Diagrammatic Language Constructs (1)}
		\label{fig:dlangspecific1}
	\end{figure}
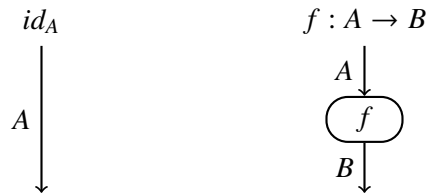

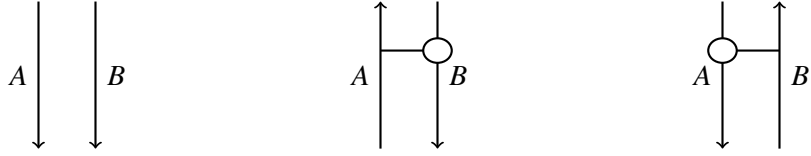
\begin{figure}[htbp]
		\centering
			\begin{tikzpicture}[thick,xscale=.75, yscale=0.65]

			\tikzset{func/.style={shape=rectangle,rounded corners=8,minimum width=1cm,draw}}
			\tikzset{claspnode/.style={shape=circle,minimum width=0.25cm,fill=white,draw}}

			\draw[->] (-6,2) -- node [left] {$A$} (-6,-1);
			\draw[->] (-5,2) -- node [right] {$B$} (-5,-1);

			% A => B

			\draw[->] (0,-1) -- node [left] {$A$} (0,2);
			\draw[->] (1,2) -- node [right] {$B$} (1,-1);
			\filldraw[fill=white] (1,1) circle (0.25);
			\draw (0,1) -- (0.75,1);

			% A <= B

			\draw[<-] (6,-1) -- node [left] {$A$} (6,2);
			\draw[<-] (7,2) -- node [right] {$B$} (7,-1);
			\filldraw[fill=white] (6,1) circle (0.25);
			\draw (6.25,1) -- (7,1);

			\end{tikzpicture}
		\caption{Basic Diagrammatic Language Constructs (2)}
		\label{fig:dlangspecific2}
	\end{figure}

	The basic constructs of the diagrammatic language of a  monoidal bi-closed category are shown in Figures \ref{fig:dlangspecific1}  and \ref{fig:dlangspecific2}. We read these diagrams  from top to bottom. Objects of the category are depicted by  arrows. For instance the left hand side arrow of Figure \ref{fig:dlangspecific1}  depicts  object $A$.  Morphisms of the category are depicted by   `blobs' with   `input'  and  `output' arrows standing for their domain and codomain objects. The right hand side blob of Figure \ref{fig:dlangspecific1} depicts  morphism $f$ which has  domain $A$ and codomain $B$. Since an identity morphism has the same domain and codomain, it is depicted in the same way as the object corresponding to it. The identity morphism $id_A$ is denoted by the same diagram as that of object $A$, e.g in the left hand side diagram of Figure \ref{fig:dlangspecific1}.

	The tensor product of two objects corresponds to the side by side placement of their arrows, for instance the tensor product  $A \otimes B$ of two objects $A$ and $B$  is depicted by putting    the arrow depicting object $A$ beside the arrow depicting object $B$, as shown in the left hand side diagram of Figure \ref{fig:dlangspecific2}.  Right and left implications  $A \multimap B$ and $A \multimapinv B$  are depicted by  side by side placement and clasping of their arrows,  pointing in opposite directions. For instance, the right implication $A \multimap B$ is depicted by  putting the arrow corresponding to    $A$ pointing upwards with the arrow corresponding to  $B$ pointing downwards; with the arrow of $B$ being clasped to the arrow of $A$. The intuition behind the directions of the arrows is that in a right implication, the antecedent has a negative role and the consequence has a positive role;  this reversal of roles is the reason behind the classical logical equivalence $A \to B \cong \neg A \vee B$. The clasp is a restriction: it is there to make us treat the implication  as one entity and to prevent us from treating each of the arrows separately. As a result, for example, we  cannot  apply   a function  to one arrow, e.g. $A$ and not the other, e.g. $B$. The diagram corresponding to the left implication $A \multimapinv B$ is the opposite of that of the right implication. These two implications are depicted in the last two diagrams of Figure \ref{fig:dlangspecific2}.

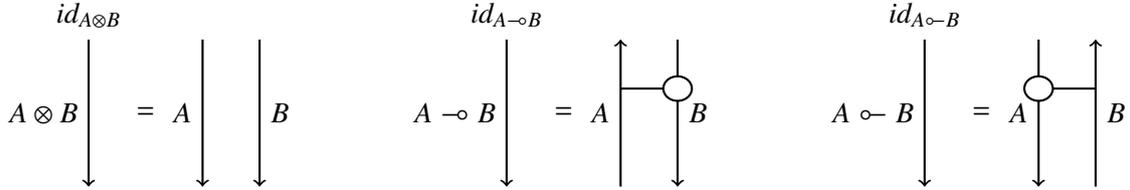
\begin{figure}[htbp]
		\centering
			\begin{tikzpicture}[thick,xscale=.75, yscale=0.65]
\hspace{-1cm}
			\tikzset{func/.style={shape=rectangle,rounded corners=8,minimum width=1cm,draw}}
			\tikzset{claspnode/.style={shape=circle,minimum width=0.25cm,fill=white,draw}}

			\draw[->] (-8,2) -- node [left] {$A \otimes B$} (-8,-1);

			\draw (-8,2.5) node {$id_{A \otimes B}$};
			\draw (-7,0.5) node {$=$};

			\draw[->] (-6,2) -- node [left] {$A$} (-6,-1);
			\draw[->] (-5,2) -- node [right] {$B$} (-5,-1);

\hspace{1cm}
			% A => B

			\draw[->] (-2,2) -- node [left] {$A \multimap B$} (-2,-1);

			\draw (-2, 2.5) node {$id_{A \multimap B}$};

			\draw (-1, 0.5) node {$=$};

			\draw[->] (0,-1) -- node [left] {$A$} (0,2);
			\draw[->] (1,2) -- node [right] {$B$} (1,-1);
			\filldraw[fill=white] (1,1) circle (0.25);
			\draw (0,1) -- (0.75,1);

\hspace{1cm}
			% A <= B

			\draw[->] (4,2) -- node [left] {$A \multimapinv B$} (4,-1);

			\draw (4, 2.5) node {$id_{A \multimapinv B}$};

			\draw (5, 0.5) node {$=$};

			\draw[<-] (6,-1) -- node [left] {$A$} (6,2);
			\draw[<-] (7,2) -- node [right] {$B$} (7,-1);
			\filldraw[fill=white] (6,1) circle (0.25);
			\draw (6.25,1) -- (7,1);

			\end{tikzpicture}
		\caption{Equivalent Diagrammatic Representations}
		\label{fig:dlangspecific3}
	\end{figure}

There might be more than one way of representing an object or morphism in this diagrammatic calculus. If so, these  different ways are considered to be `equal' and their resulting  representations  are considered to be `equivalent'. Figure \ref{fig:dlangspecific3} shows three  of the main Baez-Stay diagrammatic equivalences.  They express the fact that one can either draw the left hand side or  equivalently the right hand side diagram to represent their corresponding morphism.  For instance, one can either draw an arrow labeled with $A \otimes B$ for the tensor object $A \otimes B$ (or the identity arrow on it, i.e. $id_{A \otimes B}$) or two side-by-side arrows, one labeled with $A$ and the other with $B$.  These equivalences are not the same as the equivalences of other diagrammatic calculi \cite{SelingerSurvey}.  For instance, there is no connection between these and  the  coherence conditions of the category. 

% The two right hand side diagrams represent  part of the residuation between the tensor product and the two implications.       

\begin{figure}[h]
		\centering
			\begin{tikzpicture}[thick,scale = 0.8]

			% evall
 \hspace{-1cm}
 
			\node [func,minimum width=2cm] (evall) at (0,.25) {$ev^l_{A,B}$};

			\draw[->] (-.5,2.5) -- node[left] {$A$} (-.5,.7);
			\draw[->] (.5,2.5) -- node[right] {$A \multimap B$} (.5,.7);
			\draw[->] (evall) -- node[right] {$B$} (0,-2);

			\draw (-1,3) node {$ev^l_{A,B} : A \otimes (A \multimap B) \to B$};
			\draw (2,0.25) node {$=$};

			\draw[->] (3,2.5) -- node[left] {$A$} (3,1);
			\draw[->] (4,1) -- node[left] {$A$} (4,2.5);
			\draw[->] (5,2.5) -- node[right] {$\, B$} (5,1);
			\node[claspnode] (c1) at (5,2) {};
			\draw (4,2) -- (c1);

			\draw (3,1) .. controls +(0,-1) and +(0,-1) .. (4,1);
			\draw (5,1) .. controls +(0,-1) and +(0,1) .. (4,-.5);
			\draw[->] (4,-.5) -- (4,-2);

			\draw (4,0.25) circle (1.15cm);

 \hspace{2cm}
			% evalr

			\node [func,minimum width=2cm] (evalr) at (8,.25) {$ev^r_{A,B}$};

			\draw[->] (7.5,2.5) -- node[left] {$B \multimapinv A$} (7.5,.66);
			\draw[->] (8.5,2.5) -- node[right] {$A$} (8.5,.66);
			\draw[->] (evalr) -- node[right] {$B$} (8,-2);

			\draw (7,3) node {$ev^r_{A,B} : (B \multimapinv A) \otimes A \to B$};
			\draw (10,0.25) node {$=$};

			\draw[->] (13,2.5) -- node[right] {$A$} (13,1);
			\draw[->] (12,1) -- node[right] {$A$} (12,2.5);
			\draw[->] (11,2.5) -- node[left] {$B\ $} (11,1);
			\node[claspnode] (c2) at (11,2) {};
			\draw (12,2) -- (c2);

			\draw (13,1) .. controls +(0,-1) and +(0,-1) .. (12,1);
			\draw (11,1) .. controls +(0,-1) and +(0,1) .. (12,-.5);
			\draw[->] (12,-.5) -- (12,-2);

			\draw (12,0.25) circle (1.15cm);

			\end{tikzpicture}
		\caption{Diagrams  for Left and Right Evaluation}
		\label{fig:evrewrite}
	\end{figure}
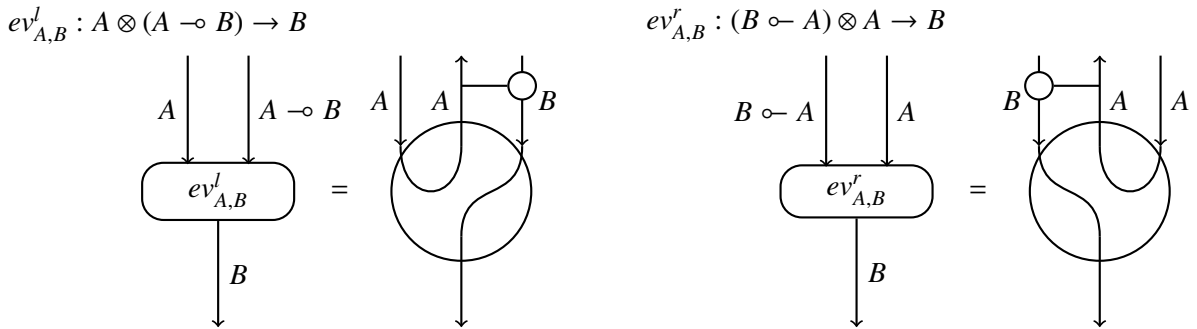

The  equivalences  for the left and right evaluations $ev^l$ and $ev^r$  are shown in  Figure \ref{fig:evrewrite}. Each of  the diagrams of  each equivalence  input two arrows and output one  arrow.  The left hand side diagrams do not depict any  information flow: they just show  which objects the evaluation morphism is being applied to. The internal structure of the evaluation morphism is being depicted inside the blob of the right hand side diagrams. These diagrams depict how  the information flows between an object  and the implication, leading to an output.   The flow represents the fact that, for instance, for the tensor product  $A \otimes (A \multimap B)$ to be evaluated and to produce the output $B$,    information has to flow from object $A$  to the antecedent of  $A \multimap B$, and similarly for the left implication. 	

\begin{figure}[h]
\begin{minipage}{6cm}
                   \begin{tikzpicture}[thick,scale = 0.7]

                   % evall

                   \node [func,minimum width=1cm] (evall) at (0,.25) {$f$};

                   \draw[->] (-.3,2.5) -- node[left] {$A$} (-.3,.7);
                   \draw[->] (.3,2.5) -- node[right] {$C$} (.3,.7);
                   \draw[->] (evall) -- node[right] {$B$} (0,-2);

                   \draw (0,3) node {$f : A \otimes C \to B$};
                   
                   \draw (5,3) node {$\Lambda^l(f) :  C \to A \multimap B$};

                   \node [func,minimum height=1.7cm,minimum width=2.2cm]  (f) at (5,.5) {\ };
                   \draw[->] (5.7,2.5) -- node {} (5.7,.55);
                               \node at (6,2.3) {$C$};
                   \node [func,minimum width=1cm] (evall) at (5.5,.1) {$f$};
                   \draw[<-] (4.23,.5) -- node[left] {$A$} (4.23,-3);
                   \draw[->] (evall) -- node[right] {$\, B$} (5.5,-3);
                   \node[claspnode] (c1) at (5.5,-2) {};
                   \draw (4.23,-2) -- (c1);
                   \draw (4.22,.45) .. controls +(0,1) and +(0,1) .. (5.2,.5);

                   \end{tikzpicture}
\end{minipage}\hspace{2cm}
\begin{minipage}{6cm} 
           \begin{tikzpicture}[thick,scale = 0.7]

                   \node [func,minimum width=1.1cm] (evall) at (0,.25) {$g$};

                   \draw[->] (-.3,2.5) -- node[left] {$C$} (-.3,.62);
                   \draw[->] (.3,2.5) -- node[right] {$B$} (.3,.62);
                   \draw[->] (evall) -- node[right] {$A$} (0,-2);

                   \draw (0,3) node {$g : C \otimes B \to A$};
        
                   \draw (5,3) node {$\Lambda^r(g) :  C \to A \lto B$};

                   \node [func,minimum height=1.7cm,minimum width=2.2cm]  (g) at (5,.5) {\ };
                   \draw[->] (4.2,2.5) -- node[left] {} (4.2,.49);
                   \node at (4.5,2.3) {$C$};
                   %\draw[->] (3.9,1.3) -- node [right] {$C$} (3.6,.5);
                   \node [func,minimum width=1cm] (evall) at (4.5,.1) {$g$};

                   \draw[->] (evall) -- node[left] {$A\ $} (4.5,-3);
                   \draw[<-] (5.62,.5) -- node[right] {$\, B$} (5.62,-3);
                   \node[claspnode] (c1) at (4.5,-2) {};
                   \draw (5.6,-2) -- (c1);
                   \draw (4.6,.49) .. controls +(0,1) and +(0,1) .. (5.62,.45);
  
        \end{tikzpicture}
\end{minipage}
\caption{Diagrams for Left and Right Currying}
\label{fig:curryeval}
\end{figure}
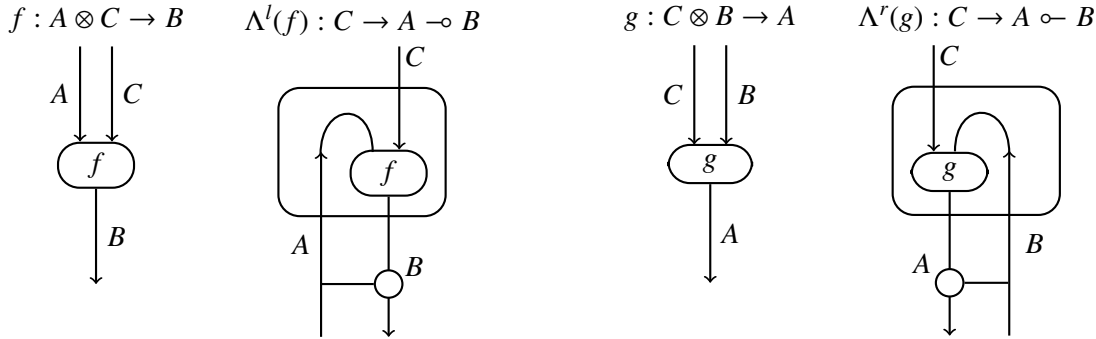

 The diagrams for left and right currying  are shown in Figure \ref{fig:curryeval}.    Here, the left hand side diagrams are the morphisms that are being curried and the right hand side diagrams depict the flow of information that happens in the process of  currying. What happens outside the blob of the right hand side diagram  has no flow and just lists the objects to which the currying morphism is being applied. The internal structure of the blobs  depict the corresponding information flows. For instance,  the diagram of $\Lambda^l(f)$ shows that in order to produce the morphism  $C \to A \multimap B$, the information encoded in $A$ and $C$ have to flow to $f\colon (A \otimes C) \to B$.  For this to happen, first $A$ has to detach itself from the clasp $A \multimap B$, a process that can only happen inside a blob. 

\begin{figure}[h] 
\centering
    \begin{tikzpicture}[thick,scale = 0.9]
    
      \begin{scope}

        \draw[->] (3,3) -- node[left] {$A$} (3,1);
        \draw[->] (4,1) -- node[left] {$A$} (4,2.5);
        \draw[->] (5,2.5) -- node[right] {$\, B$} (5,1);
        \node[claspnode] (c1) at (5,2) {};
        \draw (4,2) -- (c1);
        
        \draw (4.5,3) node (f1) [func,minimum width=.5cm] {$f$};
        \draw (4,2.5) .. controls +(0,.25) and +(-0.15,0) .. (4.22,3);
        \draw (5,2.5) .. controls +(0,.25) and +(0.15,0) .. (4.78,3);
        \draw (4.5,3) circle (.65cm);

        \draw (3,1) .. controls +(0,-1) and +(0,-1) .. (4,1);
        \draw (5,1) .. controls +(0,-1) and +(0,1) .. (4,-.5);
        \draw[->] (4,-.5) -- (4,-2);

        \draw (4,0.25) circle (1.15cm);
        
        \draw (6,1) node {$=$};
        
        \draw (7,1) node [func,minimum width=.5cm] (f1p) {$f$};
        \draw [->] (7,3) -- node [right] {$A$} (f1p);
        \draw [->] (f1p) -- node [right] {$B$} (7,-2);
        
        \draw (8,1) node {$=$};
        
        \draw[->] (11,3) -- node[right] {$A$} (11,1);
        \draw[->] (10,1) -- node[right] {$A$} (10,2.5);
        \draw[->] (9,2.5) -- node[left] {$B\ $} (9,1);
        \node[claspnode] (c2) at (9,2) {};
        \draw (10,2) -- (c2);
        
        \draw (9.5,3) node (f2) [func,minimum width=.5cm] {$f$};
        \draw (9,2.5) .. controls +(0,.25) and +(-0.15,0) .. (9.22,3);
        \draw (10,2.5) .. controls +(0,.25) and +(0.15,0) .. (9.78,3);
        \draw (9.5,3) circle (.65cm);

        \draw (11,1) .. controls +(0,-1) and +(0,-1) .. (10,1);
        \draw (9,1) .. controls +(0,-1) and +(0,1) .. (10,-.5);
        \draw[->] (10,-.5) -- (10,-2);

        \draw (10,0.25) circle (1.15cm);
      \end{scope}
      
      \begin{scope}
        
      \end{scope}
      
    \end{tikzpicture}
\caption{Diagrams   for Left and Right Yanking}
\label{fig:yank}
\end{figure}
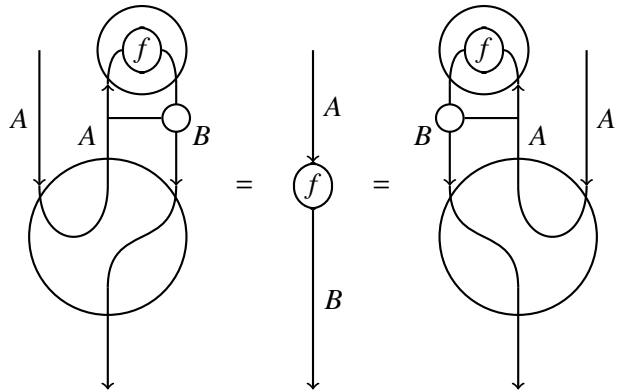

The diagrams for the equivalences corresponding to the evaluation of the names of morphisms, which lead to  yanking,  are shown in Figure \ref{fig:yank} (note that in these diagrams we have stretched  the  caps  inside the outer blobs  for reasons of space). They show how, if we first curry  $f$ to produce $A \multimap B$,  then evaluate the result with $A$, the output will be  the same as the output of  applying $f$ to $A$, that is $B$. 

\section{Monoidal Bi-Closed Diagrams and Flow of Meaning in Language}

\begin{figure}[h]
\begin{center}	
	\begin{minipage}{3cm}
		\begin{tikzpicture}[thick,scale = 0.6]
			
			\begin{scope}
				
				\draw (0, 6.5) node {Men};
				\draw (2,6.5) node {kill};

				\draw[<-] (1.5,6) -- node [left] {$n$} +(0,-1.5);
				\draw[->] (2.5,6) -- node [right] {$s$} +(0,-1.5);
				\draw (2.5,5.5) node (c4) [claspnode] {};
				\draw (1.5,5.5) to (c4);
				
				\draw[->] (0,6) -- node[left] {$n$} (0,4.5);
				
				\draw (1.5,3.5) ellipse (2cm and 1cm);
				\draw (1.5,4.5) -- (1.5,4);
				\draw (0,4.5) -- (0,4);
				\draw (0,4) .. controls +(0,-1) and +(0,-1) .. (1.5,4);
				\draw (2.5,4.5) -- (2.5,3);
				
				\draw[->] (2.5,3) -- node [right] {$s$} +(0,-1);
				
			\end{scope}			
			
		\end{tikzpicture}
\end{minipage} \hspace{3cm}
\begin{minipage}{5cm}
	\begin{tikzpicture}[thick,scale = 0.6]
		
		\begin{scope}
			
			\draw (0, 6.5) node {Men};
			\draw (2.5,6.5) node {kill};
			\draw (5,6.5) node {dogs};

			% \draw[<-] (1.5,6) -- node [left] {$N$} +(0,-1.5);
			\draw[->] (2.5,6) -- node [left] {$n \multimap s\ $} +(0,-1.5);
			\draw (2.5,5.5) node (c4) [claspnode] {};
			\draw (3.5,5.5) to (c4);
			
			\draw[<-] (3.5,6) -- node[right] {$n$} (3.5,4.5);
			\draw[->] (5,6) -- node[left] {$n$} (5,4.5);
			\draw (3.5,4) .. controls +(0,-1) and +(0,-1) .. (5,4);
			\draw (3.5,4.5) -- (3.5,4);
			\draw (5,4.5) -- (5,4);
			\draw (3.75,3.5) ellipse (2cm and 1cm);
			\draw[->] (2.5,3) -- node [left] {$n \multimap s\ $} +(0,-1);
			
			\draw[->] (0,6) -- node[left] {$n$} (0,4.5);

			\draw[dashed] (2.5,2) -- (1.5,1.5);
			\draw[dashed] (2.5,2) -- (2.5,1.5);
			
			\draw[->] (2.5,1.5) -- node [right] {$\ s$} +(0,-1);
			\draw[<-] (1.5,1.5) -- node [left] {$n$} +(0,-1);
			\draw (2.5,1) node [claspnode] (c5) {};
			\draw (1.5,1) to (c5);
			
			\draw (1.5,.5) -- (1.5,0);
			\draw (0,4.5) -- (0,0);
			\draw (0,0) .. controls +(0,-1) and +(0,-1) .. (1.5,0);
			\draw (2.5,4.5) -- (2.5,3);
			\draw (2.5,.5) -- (2.5,-.5);
			\draw (1.25,-.5) ellipse (2cm and .8cm);
			\draw[->] (2.5,-.5) -- node [right] {$s$} +(0,-1.5);

		\end{scope}			
		
	\end{tikzpicture}
	\end{minipage}
		\caption{Diagram for `men kill' and `men kill dogs'}
		\label{fig:intrans}
\end{center}
	\end{figure}
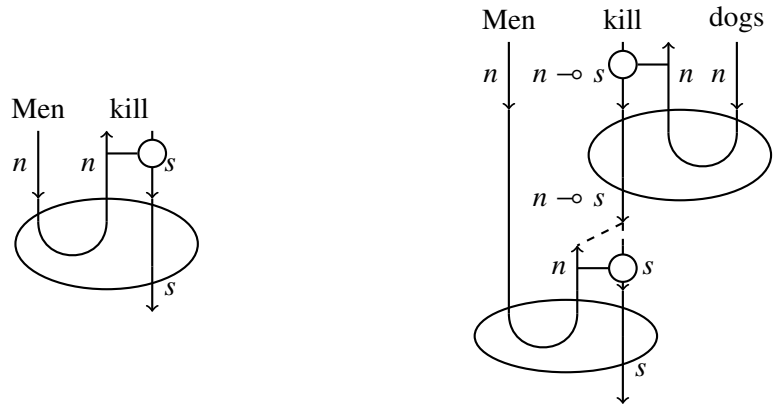

The diagram of the intransitive sentence  `men kill'  consists of  just one  evaluation: $ev^l_{n,s}$.  It shows how  information has to flow from the subject `men' to the verb `kill' to form the intransitive sentence `men kill'.   In other words, the meaning of an intransitive sentence is obtained by applying its verb to its subject. The diagram for this application is depicted in Figure~\ref{fig:intrans}. For reasons of space, we use  lowercase letters for  both the grammatical types and their interpretation.

For the transitive sentence `men kill dogs', we first depict the  subtype $n \multimap s$ of the verb with a single arrow and do an $ev^r_{n\multimap s, n}$ with the object. Then we use the equivalence for the  right implication   and replace the resulting object $n \multimap s$   with its equivalent clasped diagram  which has two arrows of types $n$ and $s$ respectively. The replacement  procedure is marked with a dotted line. At this point,  we do an  $ev^l_{n,s}$ on the implication  diagram $n \multimap s$ and the subject $n$. The full diagram in Figure~\ref{fig:intrans}, shows that to obtain a meaning for a transitive sentence,  information has to flow in two steps:  first from the  object `dogs' to the verb `kill' and then from the subject `men' to the resulting verb phrase `kill dogs'. Note that this information flow was done in one step in the compact closed setting, where the subject and object interacted with the verb simultaneously and in one step. Whereas in the monoidal case, we have to  wait for the verb to interact with the object before it can interact with the subject.

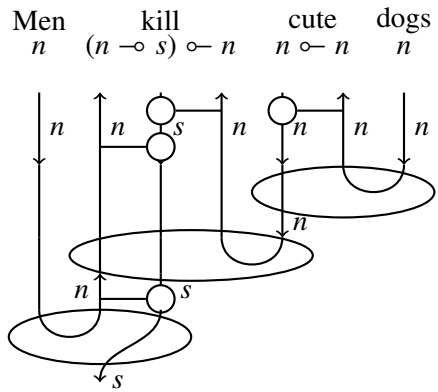
\begin{figure}[h]
\begin{center}			
			\begin{tikzpicture}
			
				\begin{scope}[yscale=0.6,thick,scale=0.8]
					%topper
					\draw (0,10) node {Men};
					\draw (2,10) node {kill};
					\draw (4.5,10) node {cute};
					\draw (6,10) node {dogs};

					\draw (0,9.25) node {$n$};
					\draw (2,9.25) node {$(n \multimap s) \lto n$};
					\draw (4.5,9.25) node {$n \lto n$};
					\draw (6,9.25) node {$n$};

%					\draw (0,8.5) node {$N$};
%					\draw (2,8.5) node {$(N \multimap S) \multimapinv N$};
%					\draw (4.5,8.5) node {$N \multimapinv N$};
%					\draw (6,8.5) node {$N$};

					%top segment

					\draw[->] (0,8) -- node[right] {$n$} ++ (0,-2);

					\draw[->] (1,6) -- node[right] {$n$} ++(0,2);
					\draw[->] (2,8) -- node[right] {$s$} ++(0,-2);
					\draw[->] (3,6) -- node[right] {$n$} ++(0,2);
					\draw (2,6.5) node (c1) [claspnode] {};
					\draw (2,7.5) node (c2) [claspnode] {};
					\draw (3,7.5) to (c2);
					\draw (1,6.5) to (c1);

					\draw[
					->] (4,8) -- node[right] {$n$} ++(0,-2);
					\draw[->] (5,6) -- node[right] {$n$} ++(0,2);
					\draw (4,7.5) node (c3) [claspnode] {};
					\draw (5,7.5) to (c3);

					\draw[->] (6,8) -- node[right] {$n$} ++(0,-2);

					%second segment

					\draw (5,6) .. controls +(0,-1) and +(0,-1) .. (6,6);
					\draw (4,6) -- (4,4.75);
					\draw[->] (4,4.75) --node[right] {$n$} (4,4);
					\draw (5,5.25) ellipse (1.5cm and .7cm);

					\draw (0,6) -- +(0,-2);
					\draw (1,6) -- +(0,-2);
					\draw (2,6) -- +(0,-2);
					\draw (3,6) -- +(0,-2);

					%third segment

					\draw (4,4) .. controls +(0,-1) and +(0,-1) .. (3,4);

					\draw (0,4) -- +(0,-2);
					\draw (1,4) -- +(0,-1);
					\draw (2,4) -- +(0,-1);

					\draw[<-] (1,3) -- node [left] {$n$} +(0,-1);
					\draw[->] (2,3) -- node [right] {$\ s$} +(0,-1);
					\draw (2,2.3) node (c4) [claspnode] {};
					\draw (1,2.3) -- (c4);

					\draw (2.5, 3.5) ellipse (2cm and .7cm);

					% last segment

					\draw (0,2) .. controls +(0,-1) and +(0,-1) .. (1,2);
					\draw[->] (2,2) .. controls +(0,-1) and +(0,1) .. (1,0) node[right] {$s$};

					\draw (1,1.25) ellipse (1.5cm and .75cm);
				\end{scope}
			\end{tikzpicture}
		\caption{Diagram for `men kill cute dogs'}
		\label{fig:adj}
\end{center}
	\end{figure}

\begin{figure}[h]
\begin{center}
\begin{tikzpicture}[yscale=0.6,thick,scale=0.8]

                               \begin{scope}[yshift=2cm]
                                  % \draw (0,10) node {Men};
                                  \draw (2,10) node {kill};
                                  \draw (4.5,10) node {cute};
                                  \draw (6,10) node {dogs};

                                  % \draw (0,9.25) node {$N$};
                                 % \draw (2,9.25) node {$(N \backslash S) / N$};
                                  %\draw (4.5,9.25) node {$N \backslash N$};
                                  %\draw (6,9.25) node {$N$};

                                  % \draw (0,8.5) node {$N$};
                                  \draw (2,8.5) node {$(n \multimap s) \multimapinv n$};
                                  \draw (4.5,8.5) node {$n \multimapinv n$};
                                  \draw (6,8.5) node {$n$};
                               \end{scope}

                               \begin{scope}

                                       %top segment

                                       % \draw[->] (0,9) -- node[right] {$N$} ++ (0,-3);
                                       %
                                       %                                       \draw[->] (1,6) -- node[right] {$N$} ++(0,3);
                                       \draw[->] (2,9) -- node[left] {$n \multimap s\ $} ++(0,-4.5);
                                       \draw[->] (3,6) -- node[right] {$n$} ++(0,3);
                                       % \draw (2.35,7.85) node {$S$};
                                       % \draw (2,6.5) node (c1) [claspnode] {};
                                       \draw (2,7.5) node (c2) [claspnode] {};
                                       \draw (3,7.5) to (c2);
                                       % \draw (1,6.5) to (c1);

                                       \draw[-] (4,9) -- (4,8);
                                       \draw[->] (4,8) -- node[right] {$n$} ++(0,-2);
                                       \draw[-] (5,6) -- node[right] {$n$} ++(0,2);
                                       \draw[->] (5,8) -- (5,9);
                                       \draw (4,7.5) node (c3) [claspnode] {};
                                       \draw (5,7.5) to (c3);

                                      % \draw (4.5,9) node (f) [func,minimum width=.5cm] {$\overline{\text{\small cute}}$};
%                                       \draw (4,8) .. controls +(0,1) and +(-.25,0) .. (4.2,9);
%                                       \draw (5,8) .. controls +(0,1) and +(.25,0) .. (4.8,9);
%                                       \draw (4.5,9) ellipse (.75cm and .9cm);

                                       \draw[->] (6,9) -- node[right] {$n$} ++(0,-3);

                                       %second segment

                                       \draw (5,6) .. controls +(0,-1) and +(0,-1) .. (6,6);
                                       \draw (4,6) -- (4,4.75);
                                       \draw[->] (4,4.75) --node[right] {$n$} (4,4);
                                       \draw (5,5.25) ellipse (1.5cm and .7cm);

                                       % \draw (0,6) -- +(0,-2);
                                       % \draw (1,6) -- +(0,-2);
                                       \draw[->] (2,4.5) -- node[left] {$n \multimap s\ $} +(0,-3.5);
                                       \draw (3,6) -- +(0,-2);

                                       %third segment

                                       \draw (4,4) .. controls +(0,-1) and +(0,-1) .. (3,4);

                                       % \draw (0,4) -- +(0,-2);
                                       % \draw (1,4) -- +(0,-1);
                                       \draw (2,4) -- +(0,-1);

                                       % \draw[<-] (1,3) -- node [left] {$N$} +(0,-1);
                                       % \draw[->] (2,3) -- node [right] {$\ S$} +(0,-1);
                                       % \draw (2,2.3) node (c4) [claspnode] {};
                                       % \draw (1,2.3) -- (c4);

                                       \draw (3, 3.5) ellipse (1.5cm and .7cm);

                                       % last segment

                                       % \draw (0,2) .. controls +(0,-1) and +(0,-1) .. (1,2);
                                       % \draw[->] (2,2) .. controls +(0,-1) and +(0,1) .. (1,0) node[right] {$S$};

                                       % \draw (1,1.25) ellipse (1.5cm and .75cm);
                               \end{scope}

                               \begin{scope}[very thick]
                                       \draw (5.55,4) -- (5.8,4) -- (5.8,.75) -- (5.55,.75);
                                       \draw[->] (6,3.5) -- +(1.75,3);

                                       \draw (8.25,9.25) -- (8,9.25) -- (8,5.75) -- (8.25,5.75);
                               \end{scope}

                               \begin{scope}[xshift=9cm]
                                          \begin{scope}[yshift=2cm]
                                               \draw (0,10) node {Men};
                                                  \draw (2,10) node {kill};
                                                  \draw (4.5,10) node {cute};
                                                  \draw (6,10) node {dogs};

                                                 % \draw (0,9.25) node {$N$};
                                                 % \draw (2,9.25) node {$(N \backslash S) / N$};
                                                 % \draw (4.5,9.25) node {$N \backslash N$};
                                                 % \draw (6,9.25) node {$N$};

                                                  \draw (0,8.5) node {$n$};
                                                  \draw (2,8.5) node {$(n \multimap s) \multimapinv n$};
                                                  \draw (4.5,8.5) node {$n \multimapinv n$};
                                                  \draw (6,8.5) node {$n$};
                                          \end{scope}

                               \begin{scope}

                                       \begin{scope}[dashed]
                                               \draw (3.5,7.5) ellipse (2.5cm and 1cm);
                                       \end{scope}

                                       \draw[dashed] (1.5,9) -- +(0,-1.5);
                                       \draw[dashed] (2.5,9) -- +(0,-1.5);

                                       \draw[dashed] (5,8) .. controls +(0,-1) and +(0,-1) .. (4,8);
                                       \draw[dashed] (5,9) -- (5,8);
                                       \draw[dashed] (4,9) -- (4,8);

                                       \draw[<-] (1.5,7.5) -- node [left] {$n$} +(0,-3);
                                       \draw[->] (2.5,7.5) -- node [right] {$s$} +(0,-3);
                                       \draw (2.5,5.5) node (c4) [claspnode] {};
                                       \draw (1.5,5.5) to (c4);

                                       \draw[->] (0,9) -- node[left] {$n$} (0,6);

                                       \draw (1.5,3.5) ellipse (2cm and 1cm);
                                       \draw (1.5,4.5) -- (1.5,4);
                                       \draw (0,6) -- (0,4);
                                       \draw (0,4) .. controls +(0,-1) and +(0,-1) .. (1.5,4);
                                       \draw (2.5,4.5) -- (2.5,3);

                                       \draw[->] (2.5,3) -- node [left] {$s$} +(0,-2);

                               \end{scope}
                               \end{scope}
                       \end{tikzpicture}
    \caption{Diagram for the Meaning of `men kill cute dogs'}
    \label{fig:adjsen}
  \end{center}
  \end{figure}
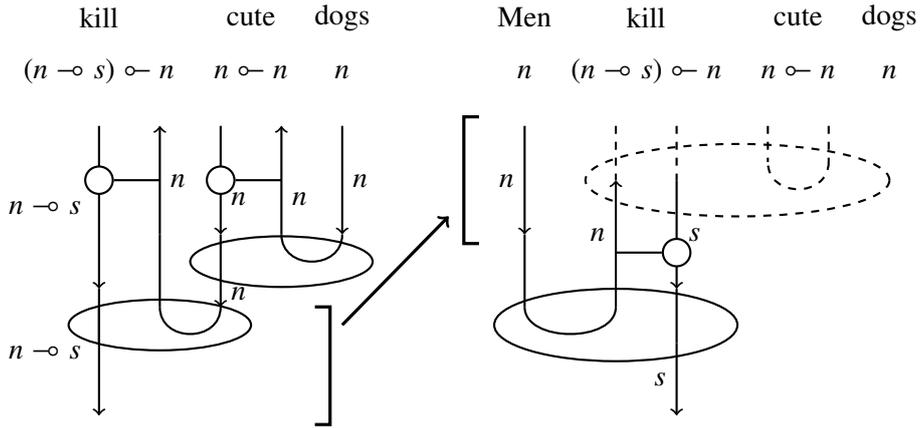

The diagram for the meaning of `men kill cute dogs' is obtained by three evaluations, as depicted in Figure~\ref{fig:adj}. Here, first the adjective acts on the object, then the verb acts on the resulting adjective-noun object to produce a verb phrase, which is then applied to the subject to produce a sentence.  The corresponding  diagrammatic computations are  depicted in Figure~\ref{fig:adjsen}. The left hand side diagram shows the verb phrase part of this interaction, whereby, `cute' is applied to `dogs' and then input to `kill' as its object. The result is then applied to the meaning of `men', this third application is depicted in the right hand side diagram, where  the dotted area stands for the general structure of the first and second  applications and serves as a  shorthand for the left hand side diagram. The arrow between the two diagrams shows that the bottom left side of the first diagram is being inserted into the right hand side and acted upon.

\begin{figure}[h]
\begin{center}
            \begin{tikzpicture}[thick,scale = 0.6]
            \draw (2,3) node {$\overrightarrow{\text{does}}\simeq$};
                   \node [func,minimum height=1.7cm,minimum width=2.2cm]  (f) at (5,.5) {\ };
                 
                   \node [func,minimum width=1cm] (evall) at (5.5,.1) {$1_{n \multimap s}$};
                   \draw[<-] (4.23,.5) -- node[left] {$n \multimap s$} (4.23,-3);
                   \draw[->] (evall) -- node[right] {$\ n \multimap s$} (5.5,-3);
                   \node[claspnode] (c1) at (5.5,-2) {};
                   \draw (4.23,-2) -- (c1);
                   \draw (4.22,.45) .. controls +(0,1) and +(0,1) .. (5.2,.6);
              
\end{tikzpicture}
 \hspace{3cm}
           \begin{tikzpicture}[thick,scale = 0.6]
            \draw (2,3) node {$\overrightarrow{\text{not}}\simeq$};
                   \node [func,minimum height=1.7cm,minimum width=2.2cm]  (f) at (5,.5) {\ };
                 
                   \node [func,minimum width=1cm] (evall) at (5.5,.1) {$\overline{\text{not}}$};
                   \draw[<-] (4.23,.5) -- node[left] {$n \multimap s$} (4.23,-3);
                   \draw[->] (evall) -- node[right] {$\ n \multimap s$} (5.5,-3);
                   \node[claspnode] (c1) at (5.5,-2) {};
                   \draw (4.23,-2) -- (c1);
                   \draw (4.22,.45) .. controls +(0,1) and +(0,1) .. (5.2,.6);

\end{tikzpicture}
\end{center}
		\caption{Diagrams for  `does' and  `not'}
		\label{fig:cutedoesnot}
	\end{figure}
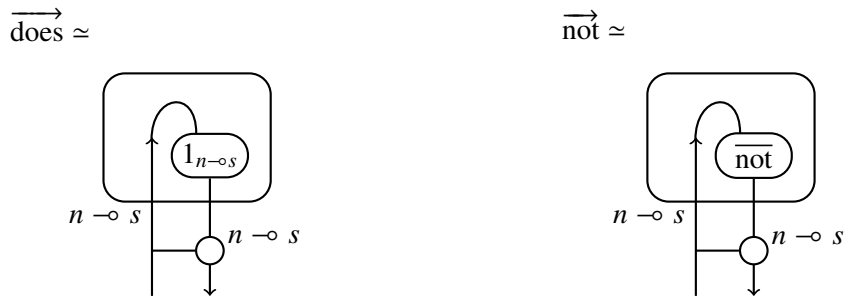

The words whose meaning vectors are \emph{names} of  linear maps are obtained by currying the corresponding linear map. For instance, the meaning vector of the auxiliary `do'   is the name of the identity and the meaning vector of the negation word `not' is the name of a negation operator,  as specified in Section~\ref{quantfunct}. The diagrams corresponding to the meaning vector of words that are names of certain morphisms are obtained by currying those morphisms. The diagrams of `does' and `not' are  depicted in Figure     \ref{fig:cutedoesnot}.

\begin{figure}[h]
\hspace{-1cm}\begin{minipage}{10cm}
  \begin{tikzpicture}[thick]
    
    \begin{scope}
      
      \begin{scope}[yshift=1cm]
        \draw (-2, 6.5) node {Men};
        \draw (.5,6.5) node {do};
        \draw (3,6.5) node {not};
        \draw (6,6.5) node {kill};
        \draw (8,6.5) node {dogs};
      \end{scope}

      % \draw[<-] (1.5,6) -- node [left] {$N$} +(0,-1.5);
      \draw[->] (6,6) -- node [left] {$n \multimap s\ $} +(0,-1.5);
      \draw (6,5.5) node (c4) [claspnode] {};
      \draw (7,5.5) to (c4);
      
      \draw[<-] (7,6) -- node[right] {$n$} (7,4.5);
      \draw[->] (8,6) -- node[left] {$n$} (8,4.5);
      \draw (7,4) .. controls +(0,-1) and +(0,-1) .. (8,4);
      \draw (7,4.5) -- (7,4);
      \draw (8,4.5) -- (8,4);
      \draw (7,3.5) ellipse (2cm and .75cm);
      \draw[->] (6,3) -- node [left] {$n \multimap s\ $} +(0,-1);
      \draw (6,4.5) -- (6,3);

      \draw (3,6.5) node [func,minimum width=.5cm] {$\overline{\text{not}}$};
      \draw (2.5,5.5) node [claspnode] (c9) {};
      \draw (3.5,5.5) to (c9);
      
      \draw (2.5,6) .. controls +(0,.25) and +(-.15,0) .. (2.62,6.5);
      \draw (3.5,5.95) .. controls +(0,.25) and +(.15,0) .. (3.38,6.5);
      
      \draw (3,6.5) circle (.6cm);
      
      \draw[->] (2.5,6) -- node[left] {$n \multimap s$} +(0,-2);
      \draw[<-] (3.5,6) -- node[right] {$n \multimap s$} +(0,-2);
      \draw (2.5,4) -- (2.5,1);
      \draw (3.5,4) -- (3.5,1.5);
      \draw (3.5,1.5) .. controls +(0,-1) and +(0,-1) .. (6,1.5);
      \draw (6,2) -- (6,1.5);
      \draw (4.25,1.25) ellipse (2cm and .6cm);
      
      \draw[->] (2.5,1) -- node [right] {$n \multimap s$} (2.5,0);

      \draw (.5,6.5) node [func,minimum width=.5cm] {$1$};
      \draw (0,5.5) node [claspnode] (c8) {};
      \draw (1,5.5) to (c8);
      
      \draw (0,6) .. controls +(0,.25) and +(-.15,0) .. (0.25,6.5);
      \draw (1,5.95) .. controls +(0,.25) and +(.15,0) .. (0.75,6.5);
      
      \draw (.5,6.5) circle (.6cm);
      
      \draw[->] (0,6) -- node[left] {$n \multimap s\,$} +(0,-1);
      \draw[<-] (1,6) -- node[right] {$n \multimap s$} +(0,-1);
      
      \draw (1,0) .. controls +(0,-1) and +(0,-1) .. (2.5,0);
      \draw [->] (1,0) -- (1,3);
      \draw (1,3) -- (1,5);
      \draw [->] (0,5) -- (0,1);
      \draw (0,1) -- (0,-1);
      
      \draw (1.25,-.5) ellipse (1.5cm and 0.5cm);
      
      \draw[->] (0,-1) -- node [right] {$n \multimap s$} (0,-1.5);

      \draw[dashed] (0,-1.5) -- (-1,-2);
      \draw[dashed] (0,-1.5) -- (0,-2);
      \draw[->] (0,-2) -- node [right] {$\ s$} +(0,-1);
      \draw[<-] (-1,-2) -- node [left] {$n$} +(0,-1);
      \draw (0,-2.5) node [claspnode] (c5) {};
      \draw (-1,-2.5) to (c5);
      \draw (-1,-3) -- (-1, -3.5);
      \draw (0,-3) -- (0,-4);
      \draw[->] (0,-4) -- node [right] {$s$} (0,-5);

      \draw[->] (-2,6) -- node[left] {$n$} (-2,4.5);
      \draw (-2,4.5) -- (-2,-3.5);
      
      \draw (-2,-3.5) .. controls +(0,-.5) and +(0,-.5) .. (-1,-3.5);
      \draw (-1,-3.75) ellipse (1.25cm and 0.5cm);

      % \draw (1.5,.5) -- (1.5,0);
      %       \draw (0,4.5) -- (0,0);
      %       \draw (0,0) .. controls +(0,-1) and +(0,-1) .. (1.5,0);
      %       
      %       \draw (2.5,.5) -- (2.5,-.5);
      %       \draw (1.25,-.5) ellipse (2cm and .8cm);
      %       \draw[->] (2.5,-.5) -- node [right] {$S$} +(0,-1.5);

    \end{scope}     
  \end{tikzpicture}
  \end{minipage}\hspace{0.5cm} $\implies$ \hspace{0.5cm}
  \begin{minipage}{5cm}
  \begin{tikzpicture}[thick]
    \draw (-2,7) node {Men};
    \draw (.5,7) node {$\overline{\text{not}}(\text{kill}(-, \text{dogs}))$};
    %\draw (2,7) node {dogs};
    
    \draw[->] (2,6) -- node [right] {$n$} (2,4.5);
    \draw (2,4.5) -- +(0,-.5);
    \draw[<-] (1,6) -- node [right] {$n$} (1,4.5);
    \draw (1,4.5) -- +(0,-.5);
    \draw[->] (0,6) -- node [left] {$n \multimap s\ $} (0,4.5);
    \draw (0,4.5) -- +(0,-1);
    \draw [->] (0,3.5) -- node [left] {$n \multimap s$} (0,2.5);
    \draw (0,5.5) node (c1) [claspnode] {};
    \draw (1,5.5) to (c1);
    \draw (1,4) .. controls +(0,-.5) and +(0,-.5) .. (2,4);
    \draw (1,4) ellipse (1.25cm and .5cm);
    
    \draw (0,2) node (notfun) [func,minimum width = .5cm] {$g$};
    \draw (0,1) node (idfun) [func,minimum width = .5cm] {$1$};
    
    \draw (0,2.5) to (notfun);
    \draw (notfun) to (idfun);
    \draw [->] (idfun) -- node [right] {$n \multimap s$} (0,0);
    
    \draw[dashed] (0,0) -- (-1,-.5);
    \draw[dashed] (0,0) -- (0,-.5);
    \draw[->] (0,-.5) -- node [right] {$\ s$} +(0,-1);
    \draw[<-] (-1,-.5) -- node [left] {$n$} +(0,-1);
    \draw (0,-1) node [claspnode] (c5) {};
    \draw (-1,-1) to (c5);
    \draw (-1,-1.5) -- (-1, -2);
    \draw (0,-1.5) -- (0,-2.5);
    \draw[->] (0,-2.5) -- node [right] {$s$} (0,-3.5);

    \draw[->] (-2,6) -- node[left] {$n$} (-2,4.5);
    \draw (-2,4.5) -- (-2,-2);
    
    \draw (-2,-2) .. controls +(0,-.5) and +(0,-.5) .. (-1,-2);
    \draw (-1,-2.25) ellipse (1.25cm and 0.5cm);
  \end{tikzpicture}
\end{minipage}
  \caption{Diagram for  Meaning of `men do not kill dogs'}
  \label{fig:mendonotkilldogs}
\end{figure}
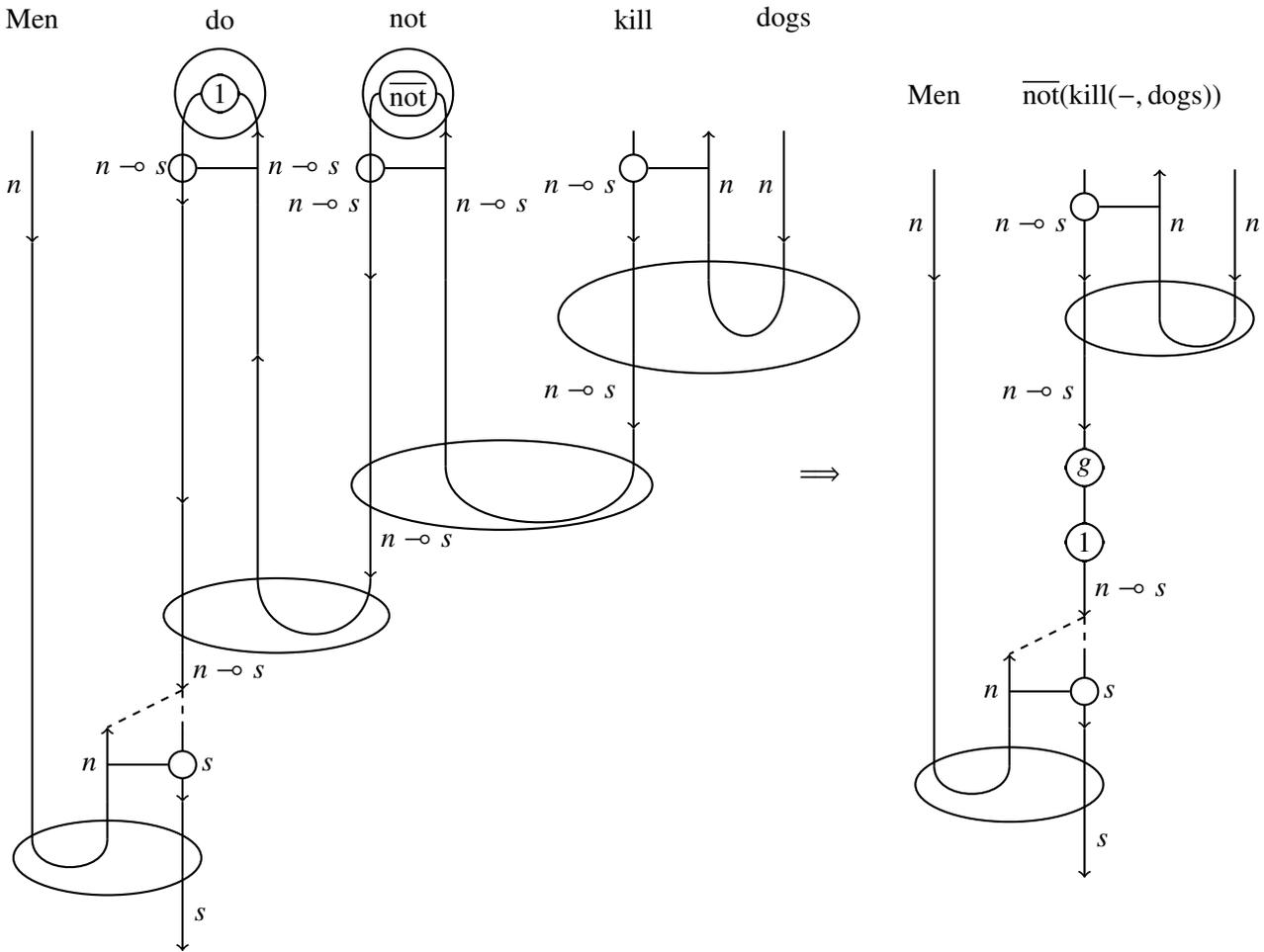

The  meaning vector  for `men do not kill dogs' is obtained by substituting the diagrams of meanings of `does' and `not' into the  diagram corresponding to the monoidal grammatical reduction  of the sentence. The result is depicted in Figure~\ref{fig:mendonotkilldogs}. Here,  the diagram on the right is the simplified version  of  the diagram on the left, where `kill'  has been applied to `dogs' and the result has been inputted to `not'. This means that the result is being negated and the negation is then being applied to the meaning of `men'. As mentioned in Section~\ref{quantfunct}, this procedure  is different from the compact closed case and this difference, which basically lies in  the order of applications,  is manifested in the corresponding diagram of each case.

If Lambek monoids are  interpreted in vector spaces, as demonstrated in the  Section~\ref{quantfunct},  they  benefit from  their extra compact structure, for instance symmetry of the tensor.  Baez-Stay diagrams can easily be turned into compact string diagrams:  by removing the  clasp restrictions and popping the bubbles. There is, however, a problem: the resulting diagrams are not always the same as the compact diagrams that were originally drawn for sentences, as demonstrated in the above example. The  problem  stems from  the restricted form of  yanking in the monoidal case.

\section{Concluding Remarks and Future Work}

We have shown that the DisCoCat framework~\cite{CCS,CSC} need  not be committed to the  grammatical formalism of  Lambek pregroups~\cite{Lambek6}, which are compact closed. Lambek monoids~\cite{Lambek}, which are monoidal closed, also allow a functorial passage to the category of finite dimensional vector spaces and the  diagrammatic calculus of Baez and Stay~\cite{Baez} can be used to depict the information flow that happens within the  sentences.  This functorial passage is used in and  referred to as   \emph{quantisation} in the context of Topological Quantum Field Theory (TQFT)~\cite{Atiyah, BaezDolan, Kock}. As future work, we would like to (1) extend the results of this paper to richer type-logics, and (2) relate their diagrammatic calculi to Baez-Stay diagrams.

With regard to (1), we would like to  develop distributional compositional models of meaning for richer grammars such as the Combinatorial Categorial Grammar~\cite{steedmanCCG} and  Lambek-Grishin algebras~\cite{Moortgat09,Bernardi}. These grammars have more expressive power and can formalise larger fragments of natural language, and in particular English; CCG has been applied to parse   large corpora of real data. To extend our quantisation functor from Lambek monoids to  CCG, the cross composition rules need to be given a categorical semantics. This would perhaps require restricting the application of these rules to subcategories of the full type-category, as in their general forms, they do not hold either in a compact or a monoidal  closed category. To extend the quantisation functor to  Lambek-Grishin algebras, we need to ask the functor to preserve the operations of the Grishin part of the algebra,  as well as their interactions with the Lambek part, hence defining functors that have to satisfy more than the monoidal conditions; the categorical properties of these functors should be studied. Whether or not such extensions prove to be trivial remains to be seen, but it is certain that augmenting the applicability of this general compositional categorical formalism will greatly assist its acceptance by the linguistic community as a practical theoretical toolkit.

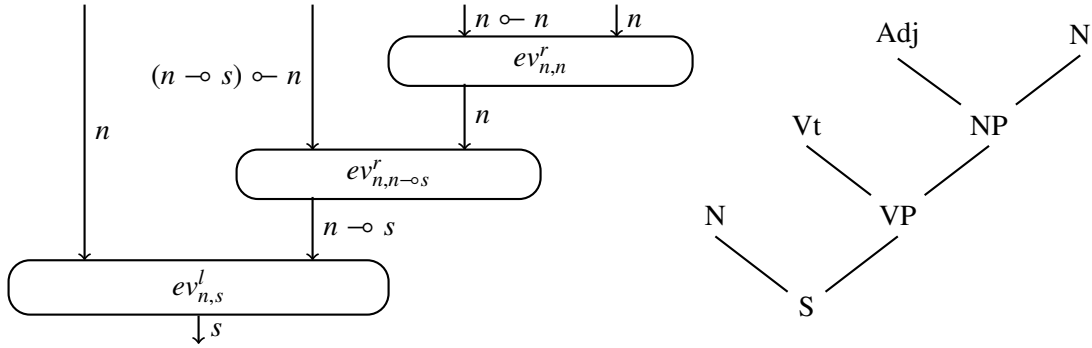
\begin{figure}[h]
		\centering
			\begin{tikzpicture}[thick]
			
			\begin{scope}[yscale=0.75]
				\draw (-.5,0) node[func,minimum width=5cm] (ev1) {$ev^l_{n,s}$}  ++(2.5,2) node[func,minimum width=4cm] (ev2) {$ev^r_{n,n \multimap s}$}  ++(2,2) node[func,minimum width=4cm] (ev3) {$ev^r_{n,n}$};

				\draw[->] (ev1) -- node[right] {$s$} (-.5,-1);

				\draw[->] (-2,5) -- node[right] {$n$} (-2,.48);

				\draw[->] (1,5) -- node[left] {$(n \multimap s) \multimapinv n$} (1,2.43);
				\draw[->] (1,1.6) -- node[right] {$n \multimap s$} (1,.48);

				\draw[->] (3,5) -- node[right] {$n \multimapinv n$} (3,4.43);
				\draw[->] (3,3.6) -- node[right] {$n$} (3,2.43);

				\draw[->] (5,5) -- node[right] {$n$} (5,4.43);
			\end{scope}
			
			\begin{scope}[xshift=7.5cm, yshift=-.25cm, scale=1.2]
				\draw (0,0) node (s1) {S} ++(-1,1) node (n1) {N} ++ (2,0) node (vp1) {VP} ++(-1,1) node (vt1) {Vt} ++(2,0) node (np1) {NP} ++(-1,1) node (adj1) {Adj} ++(2,0) node (n2) {N};
				\draw (s1) to (n1.south) ++ (s1) to (vp1.south) ++ (vp1) to (vt1.south) ++ (vp1) to (np1.south) ++ (np1) to (adj1.south) ++ (np1) to (n2.south);
			\end{scope}
			\end{tikzpicture}
		\caption{Baez-Stay Diagrams and  Parse Trees}
		\label{fig:parselikeparse}
	\end{figure}

As for (2), note that the general  structure of Baez-Stay diagrams  are the same as  parse trees. This is not surprising given the  Curry-Howard isomorphism between monoidal categories and lambda calculus~\cite{szabo}. Whereas parse trees are purely syntactical and only depict the grammatical structure of sentences, Baez-Stay diagrams have  extra information in their nodes (blobs) about the {flow of meaning} and the semantic structure of phrases. The straight lines and the  cups allow for a flow between the nodes they are connecting; this encodes  the order and details of function applications. The clasps, however,  stop this flow from  happening, hence making it explicit  which applications cannot happen. One can  identify the clasps with lines and remove the cups to obtain a  parse tree without this information. We have demonstrated this procedure for an example  sentence in  Figure~\ref{fig:parselikeparse}; here,  VP  stands for an intransitive, transitive, or ditransitive verb phrases, and  NP stands for  a noun phrase.  	

The extra information in the blobs of Baez-Stay diagrams makes them more like \emph{proof nets}. A proof net is a  diagrammatic notation  used to depict the grammatical structures   of sentences and their lambda calculus meanings. Proof nets unify proofs, hence solve the spurious ambiguity problem and   have been extensively studied  by the linear logic~\cite{Groote} and linguistics communities~\cite{Roorda, Moot, Moortgat}. They contain semantic information about sentences, represented by   \emph{traverse instructions} that describe how to form the lambda term corresponding to a grammatical reduction.  These lambda terms  can also be directly read from the sequent calculus proofs  corresponding  to the proof nets. An abstract  form of proof nets, referred to as \emph{proof structures},   encode semantic information about sentences to start with. These come equipped with certain \emph{contraction rules} and can be  rewritten  to  parse trees using these rules. The rewrites model the computations that provide us with the final meaning of the sentence.   A formal  connection between  Baez-Stay diagrams and proof nets or structures  constitutes future work.  As one of the referees suggested, one possibility would be to  develop a lambda calculus for Baez-Stay diagrams in the style of \emph{glue semantics}~\cite{dalry}.

\section{Acknowledgement}
We would like to thank Michael Moortgat, Glyn Morrill, Samson Abramsky, and Mike Stay for fruitful and clarifying discussions.

{\small 
\bibliographystyle{plain}
\bibliography{APAL-Galop} 
}
\end{document}